\documentclass{amsart}

\usepackage[T1]{fontenc}
\usepackage[utf8]{inputenc}
\usepackage{mathtools} 
\usepackage{comment}
\usepackage{tikz}
\usepackage{makecell} 

\def\hoeven#1{{#1}}

\usepackage{algorithm}
\usepackage[noend]{algpseudocode}
\algrenewcommand{\algorithmiccomment}[1]{\hfill//\hskip.1em\emph{#1}}

\usepackage[colorlinks=true,citecolor=blue]{hyperref}
\usepackage[capitalize]{cleveref} 

\newtheorem{definition}{Definition}
\newtheorem{exa}[definition]{Example}
\newtheorem{rk}[definition]{Remark}
\newtheorem{proposition}[definition]{Proposition}
\newtheorem{theorem}[definition]{Theorem}
\newtheorem{corollary}[definition]{Corollary}
\newtheorem{conjecture}[definition]{Conjecture}

\AddToHook{env/definition/begin}{\crefalias{definition}{definition}}
\AddToHook{env/exa/begin}{\crefalias{definition}{example}}
\AddToHook{env/rk/begin}{\crefalias{definition}{remark}}
\AddToHook{env/proposition/begin}{\crefalias{definition}{proposition}}
\AddToHook{env/theorem/begin}{\crefalias{definition}{theorem}}
\AddToHook{env/corollary/begin}{\crefalias{definition}{corollary}}
\AddToHook{env/conjecture/begin}{\crefalias{definition}{conjecture}}

\def\calA{{\mathcal A}}
\def\ba{{\mathbf{a}}}
\newcommand{\ord}{{\operatorname{ord}}}
\def\Lalg{L^{\mathrm{alg}}}
\def\Lmin{L^{\mathrm{min}}_f}
\newcommand\Diag{{\operatorname{Diag}}}
\def\LP{L_P}
\def\QX{{\mathbb Q}}
\def\NX{{\mathbb N}}
\def\AL{{\mathcal A}^L}
\def\Gal{\operatorname{Gal}}
\def\calP{{\mathcal P}}

\def\N{{\mathbb N}}
\def\Q{{\mathbb Q}}
\def\Z{{\mathbb Z}}
\def\Qbar{{\overline{\mathbb Q}}}
\def\Dz{{\partial_z}}

\def\pFqnoargs#1#2{{}_#1F_#2}
\def\pFq#1#2#3#4#5{\pFqnoargs{#1}{#2}\biggl[\begin{matrix}%
{#3}\\{#4}%
\end{matrix}\,;\,#5\biggr]} 

\def\sympower{{\scriptstyle \bigcirc\hspace{-0.23cm}\raisebox{0.018cm}{${\scriptstyle s}$}}\hspace{0.1cm}}

\begin{document}

\title[On deciding transcendence of power series]{On deciding transcendence of power series}

\author[A.~Bostan]{Alin~Bostan}
\address[Alin Bostan]{Inria, Sorbonne Université, LIP6, F-75252 Paris, France}
\email{alin.bostan@inria.fr}

\author[B.~Salvy]{Bruno~Salvy}
\address[Bruno~Salvy]{Inria, CNRS, ENS de Lyon, Université Claude Bernard Lyon 1,
LIP, UMR 5668, 69342 Lyon Cedex 07, France}
\email{bruno.salvy@inria.fr}

\author[M.~F.~Singer]{Michael~F.~Singer}
\address[Michael~F.~Singer]{Department of Mathematics, North Carolina State University, Box 8205, Raleigh, NC 27695, USA}
\email{singer@ncsu.edu}
 
\begin{abstract} 
	It is well known that algebraic power series are differentially finite (D-finite): they satisfy linear differential equations with polynomial coefficients. 
The converse problem---determining whether a given D-finite power series is algebraic or transcendental---is notoriously difficult.
We prove that this problem is decidable: we give two theoretical
algorithms and a transcendence test that results in an efficient
Maple implementation.
\end{abstract}

\date{\today}

\subjclass[2020]{11J81, 16S32, 34M15}


\maketitle

\section{Introduction}

\subsection{Stanley's problem and our contribution} \label{ssec:problemS}
A power series $f\in\Q[[z]]$ is called \emph{algebraic} if it is a root of some polynomial
equation {$P(z,f(z)) = 0$}, where $P \in \mathbb{Q}[x,y] \setminus \{ 0 \}$; otherwise, $f$ is called \emph{transcendental}.
A classical result, known to Abel (1827) and possibly much older, states that any algebraic power series $f\in\Q[[z]]$  is \emph{D-finite} (or, \emph{differentially finite})---that is, it satisfies a linear differential equation with coefficients in $\Q[z]$. Moreover, 
the minimal-order nontrivial homogeneous linear differential equation satisfied by~$f$ has order at most $\deg_y(P)$ and coefficients in $\Q[z]$ of degree at most $4 \deg_x(P)\deg_y(P)^2$~\cite{BoChLeSaSc07}.

Conversely, not every D-finite power series is algebraic; for
instance, $\exp(z)$ and $\log(1-z)$ are both D-finite and transcendental. Several methods are available to prove the transcendence of $\exp(z)$ and $\log(1-z)$ (see e.g., \cite{BoCaRo24}), but in general it is notoriously difficult to decide if a given D-finite function is algebraic or transcendental.  This is the topic of our article.

In his seminal article on D-finite 
functions~\cite[\S4,~(g), page~186]{Stanley80}, Richard Stanley asked the following question:

\begin{quote}
\emph{``Given the 
differential equation
\begin{equation}\label{eq:defL}
{L}(f(z))=a_r(z)f^{(r)}(z)+\dots+a_0(z)f(z)=0
\end{equation}
with polynomial coefficients $a_i \in \Q[z]$,
together with suitable (finitely many) initial conditions, satisfied by a D-finite power series $f$, give an algorithm suitable for computer implementation for deciding whether $f$ is algebraic.''} 
\end{quote}

For instance, $f(z) \coloneq \sum_{n \geq 1} z^n/n$ is a D-finite power series represented by the second-order differential equation $(z-1)f''(z) + f'(z) = 0$ and initial conditions $f(0)=0, f'(0)=1$,
and such an algorithm would have to recognize the transcendence of $f$ starting from this data alone.

Deciding transcendence of formal power series has many motivations and applications. For
instance, in number theory a first step towards proving the transcendence of a complex number
consists in proving the transcendence of some suitable power series attached to it. Other examples
come from combinatorics, where the nature of generating functions may reveal strong
underlying structures~\cite{BoChHoKaPe17}, and from computer science, where a natural
question is whether algebraic power series are easier to manipulate than transcendental
ones~\cite{BoBo88}.

Although Stanley insisted on the practical aspect of the targeted algorithm for deciding
algebraicity of D-finite functions, not even a theoretical one has appeared so far. In this
work, we make progress on Stanley's question: we give two algorithms that prove the
decidability of the algebraic or transcendental nature of D-finite power series in $\Q[[z]]$
and we give an (incomplete) transcendence test whose implementation is efficient in practice.

The first theoretical algorithm (\cref{algo:Stanley1}) relies on minimization of linear
differential equations~\cite{BoRiSa24}; the second one 
(\cref{algo:Stanley3}) uses factorization.
Both rely on Singer's algorithm~\cite{Singer79}. The practical transcendence test
(\cref{algo:Stanley2a}) combines the same minimization as the first algorithm 
with
a local analysis; it is incomplete in the sense that it is guaranteed to be correct only when
it returns ``transcendental''. If the input differential equation possesses additional
properties (e.g., if it annihilates the diagonal of a multivariate rational function~\cite{Christol15}), then a
variant of it (\cref{algo:Stanley2b}) is also guaranteed to be correct when it returns
``algebraic'', modulo a conjecture in the theory of arithmetic differential equations
(\cref{conj:ChristolAndre}). Designing complete and efficient algebraicity or transcendence
tests is still open.

\subsection{Difficulties and analogies} \label{ssec:analogies}

There are several \emph{a priori} reasons why Stanley's problem is 
difficult.
First, the minimal polynomial 
of an algebraic power series $f(z)$
may have arbitrarily large size
(degrees) with respect to the size (order/degree) 
of its defining
differential equation. For instance, 
$f(z) = \sqrt[N]{1-z}$ has algebraicity degree~$N$
and satisfies the first-order differential equation
$N(z-1) f'(z)-f(z)=0$. 
Second, no characterization 
is available
for coefficient sequences
of algebraic power series; 
this is in contrast with the smaller class 
of rational functions
(whose coefficient sequences are C-recursive, i.e., satisfy linear recursions with constant coefficients) and with the larger class of D-finite functions
(whose coefficient sequences are P-recursive, i.e., satisfy linear recursions with polynomial coefficients)\footnote{Note that for the class of {diagonals of rational functions}, which is 
intermediate between algebraic and D-finite power series,
coefficient sequences 
correspond (conjecturally)
to the class of
P-recursive sequences with (almost) integer coefficients and geometric growth (this is predicted by the ``Christol conjecture'', see Part (1) of~\cref{conj:ChristolAndre}).}.

There are several analogies between the problems of
transcendence in $\mathbb{Q}[[z]]$
and
of irreducibility in $\mathbb{Q}[z]$.
One is that ``generic'' power series in $\mathbb{Q}[[z]]$
are transcendental,
much like 
``generic'' polynomials in $\mathbb{Q}[z]$
are irreducible.
Another one is that, 
while many
sufficient criteria exist
(e.g., {Eisenstein's} criteria for polynomial irreducibility~\cite{Eisenstein1850}, see also~\cite{Schoenemann1846,Schoenemann1850,Cox11}, and for power series transcendence, see \cref{ssec:Eisenstein}), none 
of them 
is also  necessary.
This is closely related to the fact that there are no known characterizations at the level of the coefficient sequence for recognizing either irreducibility in~$\Q[z]$, or transcendence in~$\Q[[z]]$.
However, since polynomial irreducibility is known to be {decidable}~\cite{Kronecker1883,LLL82}, 
it is legitimate to hope, by analogy, that the same holds for transcendence of power series.

\subsection{Related problems}\label{ssec:related}
Although Stanley's question is quite recent, it is related to many classical questions and results, dating back to the beginning of the 19th century.

Given a linear differential operator $L=\sum_{i=0}^r a_i(z) \partial_z^i$ 
(where $\partial_z$ stands for $\frac{d}{dz}$)
with coefficients in $\mathbb Q[z]$,
modeling the linear differential equation~\eqref{eq:defL},
one may ask several questions concerning the algebraic
nature over $\mathbb Q(z)$ of its solutions:

\begin{enumerate}
\item[{\bf (F)}] Decide if \emph{all} solutions of $L$ are algebraic (Fuchs' problem).
\item [{\bf (L)}] Decide if $L$ admits \emph{at least one} nontrivial algebraic solution (Liouville's problem).
\item [{\bf (S)}] Decide if \emph{a given} nonzero solution $f$ of $L$ is algebraic (Stanley's problem).
\end{enumerate}

\begin{rk}
In problems {\bf (F)} and {\bf (L)} an ``algebraic solution'' of~$L$ is an
element $f$ of the field of formal Puiseux series $\mathcal P=\bigcup_
{n\in
\NX^\star} \Qbar
((z^{\frac{1}{n}}))$, 
such that $f$ is algebraic over $\mathbb Q(z)$ and $L(f) = 0$. 
(Here $\Qbar$ is
the algebraic closure of $\QX$.)

``All
solutions of $L$ are algebraic'' means that $L$ admits a \emph{full
basis}, i.e., a
linearly independent family of
cardinality $r=\ord(L)$, of solutions in~$\mathcal P$.
\end{rk}

Note also that if $L$ is irreducible\footnote{This means that $L$
cannot be written $L=L_1 L_2$ in the (non-commutative) ring $\mathbb{Q}(z)\langle \partial_z \rangle$ of linear differential operators, with $1 \leq \ord(L_2) < r= \ord(L)$.
Irreducibility of linear differential operators was first studied by Frobenius (1873), and then uniqueness properties for factorization into irreducible operators were
further studied by Landau (1902) and Loewy (1903).
For historical accounts on the rich structural analogy between polynomials and linear differential equations, we refer to 
\cite{Demidov83},
\cite[\S2.4]{Singer96}
and
\cite[p.~100]{BoChLiSa12},
and the references therein.
}, then the three problems {\bf (F)}, {\bf (L)} and {\bf (S)} are equivalent.
This is because an irreducible operator either admits a full basis of
algebraic solutions, or no nontrivial algebraic solution, see
\cite[Prop.~2.5]{Singer79}.

Above, the word ``decide'' may have several meanings.
One may for instance ask for a criterion, that is for a mathematical characterization 
based on the inspection of the coefficients of~$L$, or of some
invariants of~$L$ (such as its singularities,
exponents at singularities, etc.\footnote{Finite singularities (or
singular points) of \cref{eq:defL} are points where one of the
rational functions~$a_i(z)/a_r(z)$ has a pole. Infinity is also a
singular point when~0 is a singularity of the equation
obtained by
making the change of variables~$z\mapsto1/z$. For more
information, see~\cite [\S15.1] {Ince1956}.}).
The classical meaning in computer science---the one that we adopt in
this article---is that one asks for an algorithmic decision
procedure, which, on every specific instance of~$L$, is able to answer
the given problem using a finite number of operations  
($+$, $\times$, $\div$, $=$) in~$\Q$.

In addition to such a decision procedure, one may also want to exhibit an
annihilating polynomial (usually, the minimal polynomial) of the
algebraic solution(s). 
Computing such a polynomial yields the most obvious ``algebraicity witness'',
on which the decision procedure itself may rely (see~\cref{ssec:GnP}).
However, these witnesses may have prohibitively large sizes;
deciding algebraicity of solutions does not necessarily need them.

\medskip 
In this article, we add one more problem to the previous list:
\begin{enumerate}
\item[{\bf (P)}] \label{pb-P} Compute the (monic) right 
factor\footnote{This is a linear differential operator with coefficients
in~$\mathbb Q(z)$, see \cref{prop:basis_alg}.} $\Lalg$ of
$L$ whose solution space is spanned by the algebraic
solutions of~$L$.
\end{enumerate}
Again, if $L$ is irreducible, then {\bf (P)} is equivalent to the problems {\bf (F)}, {\bf (L)} and {\bf (S)}.

Note that a solution to {\bf (P)} implies solutions to each of the problems {\bf (F)}, {\bf (L)} and~{\bf (S)}.
Indeed, solving {\bf (F)},  {\bf (L)}, and {\bf (S)}
amounts, respectively, to checking  
$\ord(L)=\ord(\Lalg)$, resp. $\Lalg \neq 1$, 
or  $\Lalg(f) = 0$
(the last equality can be checked as explained in~\cite[Lemma 2.1]{BoRiSa24}).

\subsection{Sufficient transcendence criteria}\label{ssec:criteria}
There exist several criteria that can be used to prove that a given power series $f\in\Q[[z]]$ is transcendental. They are simply built on 
properties of algebraic power series. We recall below two important ones,
Eisenstein's arithmetic criterion and Flajolet's analytic criterion.
Although very useful in practice, none of these sufficient criteria is also necessary,
and there are concrete examples that escape them 
(see \cref{ex:trident}). 
This observation brings additional motivation to Stanley's problem~{\bf (S)}.

\subsubsection{Eisenstein's criterion} \label{ssec:Eisenstein}
\begin{definition}\label{def:gb}
A power series $f = \sum_{n\geq 0} a_n z^n$ in $\Q[[z]]$ is called \emph{globally bounded}
if it has a non-zero and finite radius of convergence, and if there exists
$C\in \mathbb{N}^*$ such that  $a_nC^n\in\mathbb{Z}$ for all $n\geq 1$. 
\end{definition}	
A famous criterion, stated (and sketchily proved) in 1852
by Gotthold Eisenstein~\cite{Eisenstein1852} and fully proved one year later by Eduard Heine~\cite{Heine1853}, asserts that if
a power series $f \in \Q[[z]]$ is algebraic, then only a finite number of prime numbers can divide the denominators of the coefficient sequence $(a_n)_{n \geq 0}$.
This is also a consequence of a theorem\footnote{Heine's 1854 result is what is commonly called ``Eisenstein's criterion'' in the literature. Hence, strictly speaking there are two distinct criteria, one due to Eisenstein and the other one due to Heine.
It might well be that they coincide for D-finite functions; this is
related to Siegel's conjecture that the classes of ``broad $G$-functions'' and ``strict $G$-functions'' coincide, see e.g., \cite[Conj.~2.3.1]{Andre00}.}
by Heine~\cite{Heine1854}
stating that any algebraic power series $f \in \Q[[z]]$
is globally bounded.
This arithmetic criterion immediately implies that the power series
$\log(1-z) = -\sum_{n\geq 1} z^n/n$
(more generally, the polylogarithm $\textrm{Li}_s(z) = \sum_{n\geq 1} z^n/n^s$ for $s\in\mathbb N^\star$),
and 
$\exp(z) = \sum_{n\geq 0} z^n/n!$ are transcendental.
However, when $f(z)$ is not given in closed form, but rather like in
Stanley's problem~{\bf (S)}, by a linear differential equation with
sufficiently many initial conditions, then Eisenstein's criterion is
difficult to apply. The main reason for this is a quite fundamental one: it is currently not known (and it is considered to be a difficult open problem) how to recognize that a P-recursive sequence $(a_n)_n$ given by a linear recurrence and initial conditions has integer terms (or, that it has \emph{almost integer} coefficients, in the sense that $a_nC^n$ are integers for all $n\geq 1$, for some constant~$C\in\mathbb{N}^*$). 
The same difficulty arises in a different context for the mere problem of representing $E$-functions, see~\cite{BoRiSa24}.

\subsubsection{Flajolet's criterion} \label{ssec:Flajolet}
In many situations, 
the given D-finite power series $f(z)$ already has integer coefficients, therefore Eisenstein's criterion is useless.
This is systematically the case in combinatorics when $f(z)$ is the generating function of a sequence whose $n$th term counts the number of objects of size~$n$ in some combinatorial class. 
Much more useful is \emph{Flajolet's criterion}~\cite[Criterion D]{Flajolet87}
(see also~\cite[VII.7.2]{FlSe09}) based on the Newton-Puiseux theorem
and on Darboux's transfer results from the local behavior of $f(z)$
around its complex singularities to the asymptotic behavior of its
coefficient sequence $(a_n)_{n \geq 0}$.

One form of the criterion asserts that \emph{if
$a_n \sim \gamma \,  \beta^n \, n^r$ as $n\rightarrow\infty$, with
either 
$r$ an irrational number or a negative integer,
or
$\beta$ a transcendental number,
or
$\gamma \cdot \Gamma(r+1)$ a transcendental number,
then $f(z)$ is transcendental.}

For more general and more refined versions on the asymptotic behavior of the coefficients of algebraic power series, we refer to  \cite[p.~501]{FlSe09}, \cite{Orlov94} and \cite[Theorem~3.1]{RiRo14}.

This criterion is the most commonly used for transcendence
questions in combinatorics.
It allows for instance to show
that if $k \in \mathbb{N}^\star$, then the hypergeometric power series $\sum_{n \geq 0} \binom{2n}{n}^k z^n$
(mentioned at the end of problem~(g) in~\cite[\S4]{Stanley80}) 
is algebraic if and only if $k=1$ (see~\cref{ssec:hyp} for a different
proof), 
and also that the non-hypergeometric power series
$\sum_{n\geq0} \sum_{k=0}^n \binom{n}{k}^2 \binom{n+k}{k}^2 z^n$ 
(occurring in Apéry's celebrated proof of the irrationality of $\zeta(3)$, 
see~\cref{ex:Apery})
is transcendental. A more general example that can be handled by
Flajolet's criterion is 
\[
\sum_{n \geq 0} \left( 
\sum_{k=0}^n 
\binom{n}{k}^{p_0} 
\binom{n+k}{k}^{p_1}
\binom{n+2k}{k}^{p_2}
\dotsm
\binom{n+mk}{k}^{p_m} 
\right) z^n,
\]
with \(p_i\in\mathbb N\), that is discussed in Appendix~\ref{appendix:Apery}.

However, there are cases of D-finite power series $f(z)$ 
for which both Eisenstein's criterion and Flajolet's criterion fail to detect the transcendence of $f(z)$, see~\cref{ex:trident}. 

\subsection{A bit of history} \label{ssec:history}
Since the beginning of the 19th century, 
many works have been dedicated to solving problems {\bf (F)} and {\bf 
(L)},
using various tools, from geometry (Schwarz, Klein), to invariants (Fuchs,
Gordan) and from group theory (Jordan) to differential Galois theory (Kolchin,
Singer). 

Already in 1833, Liouville~\cite{Liouville1833} proposed an 
algorithm for computing a basis of \emph{rational} solutions (that is, solutions in $\mathbb{Q}(z)$) 
of linear differential equations like~\eqref{eq:defL}.
This algorithm (with some enhancements and improvements, e.g., 
\cite{AbKv91}) 
is the one currently implemented in most computer algebra systems.
Liouville's algorithm 
clearly solves the variants 
{\bf (F)}$_{\text{rat}}$ 
and 
{\bf (L)}$_{\text{rat}}$ 
of problems
{\bf (F)} and {\bf (L)},
in which the word ``algebraic'' is replaced by ``rational'';
it can also be used to solve 
the variant~{\bf (S)}$_{\text{rat}}$ of
problem~{\bf (S)}:
from a basis of  
solutions $r_1(z), \ldots, r_s(z)$ in~$\Q(z)$ of~$L$, one first computes 
the right factor $L^\textrm{rat}$ of $L$ whose solution space is
spanned by the rational solutions of $L$,
as the LCLM\footnote{We denote by $\text{LCLM}\{L_1, \ldots, L_k\}$
the \emph{least common left multiple} of the operators $L_1, \ldots, L_k$, that is the monic minimal order
operator whose solution space contains the solution spaces of all $L_i$'s.} 
of the differential operators $\partial_z - r_j'(z)/r_j(z)$,
and then one checks whether $L^\textrm{rat}(f)$ is zero or not (using~\cite[Lemma~2.1]{BoRiSa24}).

In 1839 Liouville addressed in~\cite{Liouville1839} what we call problem {\bf (L)}.
He partially solved it for second-order differential operators.
For instance, 
he showed that the equation $y''+ r(z)y = 0$ does not admit any nontrivial algebraic solution when $r(z) \in\Q[z] \setminus \{ 0 \}$. 
(In particular, this holds for Airy's equation $y'' = zy$.)
One year later he proved the same for $r(z)=(1+z^2)^2/(2z-2z^3)^2$, from which he deduced that the complete elliptic integral $f(z) = \int_{0}^{1} dx/\sqrt{(1-x^2)(1-z^2x^2)}$ is not algebraic 
(and not even Liouvillian, that is, solvable in terms of integrals,
exponentials and algebraic functions).
In 1841, he
applied his method to the Bessel equation
$z^2 f''(z)+zf'(z) + (z^2-\mu^2) f(z)=0$.
However, he could not solve the problem when $r(z)$ is an arbitrary rational function;
he reduced it to upper bounding the possible algebraicity degree of a solution. With such a bound $B$ at hand, Liouville reduced problem {\bf (L)} to finding rational solutions of the symmetric powers $L^{\sympower m}$ 
(that is the monic minimal order operator whose solution space
contains the powers~$f^m$ for all solutions $f$ of~$L$) for $1\leq m
\leq B$.

Liouville's work was taken up by Pépin (1863, 1878, 1881), 
who focused on problem~{\bf (F)} and 
managed to remove the restriction on the algebraicity degree~\cite{Pepin1863,Pepin1878}.
His 1863 paper contained a few errors corrected in 1878 (after remarks by Fuchs).
Pépin completed his study in his long memoir published in 1881~\cite{Pepin1881}; there,
he proved that if the equation $y''+ r(z)y = 0$ has algebraic solutions only, 
then it admits a basis of solutions $\{ y_1, y_2 \}$ such that either
$(i)$  $y_1= \sqrt[m]{a}$ and $y_2 = b/y_1$, for $a, b$ in $\Q(z)$ and $m\in\N^\star$;
or 
$(ii)$ $y_i^m$ are both roots of a quadratic equation over $\Q(z)$ for $m\in\N^\star$;
or 
$(iii)$ $y_i^m$ are both roots of an equation of degree $\mu$ over $\Q(z)$,
where $(m,\mu)\in \{ (4, 6), (6, 8), (12, 10)\}$.
In all cases, $y_i'/y_i$ is algebraic of degree 1, 2, 4, 6 or 12.
        
Meanwhile, in 1873, Schwarz~\cite{Schwarz1872} famously solved problem {\bf (F)} for
second-order operators with three singular points (the Gauss
hypergeometric equation), see~\cref{ssec:hyp}.
Fuchs showed in 1876 that if $y''+ r(z)y = 0$ has a nontrivial
algebraic solution, then there exists a binary homogeneous form $F$
of degree $d \leq 12$ such that, for some basis $\{ y_1, y_2
\}$ of solutions, $F(y_1, y_2)$ is the $k$-th root of a rational function for some explicit
number~$k$ depending only on~$d$.
Inspired by Fuchs's work, Klein showed in~1877 that any second-order
linear differential equation with only algebraic solutions comes from
some hypergeometric equation from Schwarz's list via a rational change
of variables.
In~\cite{BaDw79}, Dwork and Baldassarri discuss Klein's and Fuchs's
articles and make Klein's approach algorithmic.
Sanabria Malag\'{o}n  generalizes this approach in \cite{Sanabria22},
where he also gives a history of the evolution of this approach
through his work and that of others. Another approach was proposed by
Jordan~(1878) who showed that if all solutions of~\eqref{eq:defL} are
algebraic, then there exists a solution whose logarithmic derivative
is algebraic of degree at most~$J(r)$, an explicit number depending
only on the order~$r$. Jordan's article and the works by Painlev\'{e} 
(1887) and his student Boulanger (1898) are the starting points of
modern algorithms by Singer~\cite{Singer79,Singer81}, who completely solved problem~{\bf (F)}, see~\cref{ssec:SingerF}. 
References to improvements  on this algorithm are given in~\cite{Singer99}.

For more historical aspects related to these rich problems,
we refer the interested reader to the following sources: 
\cite[p.~1--13]{Boulanger1898},
\cite[p.~160--165]{Vessiot1910},
\cite[p.~5--10]{Gray84},
\cite[p.~407--413]{Lutzen90},
\cite[p.~25]{Singer90},
\cite[p.~533--536]{Singer99},
\cite[p.~124]{PuSi03}
and
\cite[p.~48--50, 273--274, and Chap. III]{Gray08}.

\subsection{The hypergeometric case} \label{ssec:hyp}
A very particular, yet quite important special case of D-finite power
series $f\in\Q[[z]]$ is the \emph{hypergeometric} class. This means that the coefficient sequence of $f(z)$  satisfies a first-order homogeneous linear recurrence with polynomial coefficients. Typical examples are
$\log (1-z)$, ${\arcsin(\sqrt{z})}/{\sqrt{z}}$, ${(1-z)^\alpha}$ for $\alpha\in\Q$, and more generally the Gaussian hypergeometric function with parameters $a, b, c \in\Q$
($-c\notin \N$),
\[
{\pFq21{a\ b}c z \coloneq \sum_{n=0}^\infty \frac{(a)_n(b)_n}{(c)_n} \,
\frac{z^n} {n!}, 
\quad
\text{where} \quad 
(d)_n \coloneq d(d+1)\cdots(d+n-1),
}
\]
solution of the differential equation 
$z(1-z)f''(z) +
(c-(a  +  b  +  1)z)f'(z) - abf(z) =0$.

Deciding the algebraicity of $_2F_1$ functions is an old problem, first solved by Schwarz~\cite{Schwarz1872} using geometric and complex analytic tools, 
and later by Landau~\cite{Landau1904,Landau1911} and Errera~\cite{Errera1913} using number theoretic tools.
Both approaches are algorithmic:  Schwarz's criterion reduces the problem to a table look-up after some preprocessing on the parameters $a,b,c$; the Landau-Errera criterion amounts to checking a finite number of inequalities. 
The reader is referred to~\cite[\S2.1.3]{BoCaRo24} 
and to Matsuda's book~\cite[Chap.~I, II and IV]{Matsuda85} for more details.

The definition of the Gaussian hypergeometric  function extends to the
\emph{generalized
hypergeometric function} with complex parameters $a_i$ and~$b_i$ with
$-b_i\not\in\mathbb N$
\[
{_kF_{\ell}}\biggl[\begin{matrix}
  {a_1}\kern.707em {a_2}\kern.707em{\cdots}\kern.707em{a_k}\\ 
  {b_1}\kern.707em{\cdots}\kern.707em{b_{\ell}}
\end{matrix}\,;\, z\biggr]
\coloneq \sum_{n=0}^\infty \frac{(a_1)_n\cdots (a_k)_n}{(b_1)_n\cdots 
(b_{\ell})_n} \,
\frac{z^n} {n!}. 
\]

When $\ell\neq k-1$, a $
{}_kF_\ell$ non-polynomial hypergeometric function in~$\mathbb
Q[[z]]$ cannot be
algebraic, by Stirling's formula combined with Flajolet's criterion;
alternatively, because for $k>\ell+1$ the power series ${}_kF_\ell$
has radius of convergence 0, while ${}_kF_\ell$ is an
entire function for $k\leq\ell$.

When~$\ell=k-1$ and the parameters are rational numbers, the
Landau-Errera criterion
was extended by Christol~\cite [\S VI]
{Christol86} under a conjecture, 
and by Beukers and Heckman~\cite{BeHe89} unconditionally.
To state it, we need the following definitions: 
for~$x\in\mathbb{R}$ we denote by $\langle x \rangle$ its fractional part $x - \lfloor x \rfloor$ if $x\notin \Z$, and 1 if $x\in\Z$; 
we say that two equinumerous disjoint multisets of real numbers $\{ u_1, \ldots, u_k\}$  and $\{ v_1, \ldots, v_k\}$
with $\langle u_1\rangle \leq \cdots \leq \langle u_k\rangle$
and 
$\langle v_1\rangle \leq \cdots \leq  \langle v_k\rangle$
\emph{interlace} if 
$\langle u_1\rangle  < \langle v_1\rangle < \cdots < \langle u_k\rangle  < \langle v_k\rangle$.
Now, let ${\bf a} = \{ a_1,\ldots, a_k \} \subset \Q$ and ${\bf b} = \{ b_1,\ldots, b_{k-1}, b_k=1 \} \subset \Q \setminus (-\N)$ 
be two multisets of rational parameters assumed to be disjoint modulo~$\mathbb{Z}$.
This assumption is equivalent to the irreducibility of the generalized hypergeometric operator
$H_{\bf b}^{\bf a} \coloneq  (z\partial_z+b_1-1)\cdots (z\partial_z+b_{k-1}-1) z\partial_z - z (z\partial_z+a_1)\cdots (z\partial_z+a_k)$.
Let $D$ be the common denominator of $a_1, \ldots, a_k, b_1, \ldots, b_{k-1}$. 
Then, the function ${}_kF_{k-1}$
is a solution of $H_{\bf b}^{\bf a}$ and the \emph{interlacing criterion} says that it
is algebraic if and only if  for all $1\leq \ell < D$ with $\gcd(\ell, D) = 1$ 
the multisets $\ell {\bf a}$ and~$\ell {\bf b}$ interlace 
(pictorially, this means that the points $\{ e^{2 \pi i \ell a_j}, j \leq k \}$ and $\{ e^{2 \pi i \ell b_j}, j \leq k \}$
interlace on the unit circle).

For instance, the Beukers-Heckman criterion immediately implies that the power series
\[\sum_{n=0}^\infty \binom{2n}{n}^{\!k} z^n = 
{_kF_{k-1}}\biggl[\begin{matrix}
  {\frac12}\kern.707em {\frac12}\kern.707em{\cdots}\kern.707em{\frac12}\\ 
  {1}\kern.707em{\cdots}\kern.707em{1}
\end{matrix}\,;\, 4^k \, z\biggr]
\]
is transcendental for all $k\geq 2$, since the interlacing condition is violated. 
This example is mentioned by Stanley in~\cite[\S4,~(g)]{Stanley80} as one of the motivations of his general question {\bf (S)}.

The Beukers-Heckman interlacing criterion does not cover all
cases of hypergeometric functions in~$\mathbb Q[[z]]$; it cannot be
used if some parameters are irrational nor if a difference between a top and a bottom parameter is an integer.
For this reason, the following examples escape it
\begin{equation}\label{ex:FuYu24}
\pFq21{1\ 1}2 z, \quad
\pFq21{2\ 2}1 z \quad \text{and} \quad
{_3F_{2}}\biggl[\begin{matrix}
  {1/2}\kern.707em {1+\sqrt{2}}\kern.707em{1-\sqrt{2}}\\ 
  {\sqrt{2}}\kern.707em{-\sqrt{2}}
\end{matrix}\,;\,  z\biggr] .
\end{equation}
Fürnsinn and Yurkevich~\cite{FuYu24} provided a complete classification of the algebraic generalized hypergeometric functions with no restriction on the set of their parameters, thus answering completely Stanley's question {\bf (S)} for power series whose coefficient sequences satisfy a recurrence of order~1.
Their result is algorithmic and makes it possible to prove that, among
the three examples in~\cref{ex:FuYu24},
the first power series is transcendental, while the last two are algebraic.
Their algorithm relies on an ingenious, albeit elementary, reduction of the
general case to the interlacing criterion of Beukers and Heckman~\cite[Fig.~1]{FuYu24}.

\subsection{Singer's algorithm for  problem {\bf (F)}}    \label{ssec:SingerF}
In the 1980s, Michael Singer~\cite{Singer79,Singer81} designed algorithms that solve Fuchs' problem {\bf (F)}.
In fact, the problem had essentially been solved in the second half of the 19th century, for order 2 equations by Schwarz, Klein, Fuchs, and for order 3 equations by Painlevé and Boulanger.  For an irreducible operator $L$ of order~$r$, Singer's algorithm in~\cite{Singer79} relies on differential Galois theory and proceeds in two main steps: 
\begin{enumerate}
\item[(i)] (Jordan's bound) decide if the  (nonlinear) Riccati differential equation\footnote{By definition, $R_L$ has order $r-1$ and has the property that $f$ is a nonzero solution of $L$ if and only if $f'/f$ is a solution of $R_L$.} $R_L(u)=0$ attached to~$L$ admits an algebraic solution~$u$ of degree at most $(49r)^{r^2}$; 
\item[(ii)] (Abel's problem) given an algebraic~$u$, decide whether $y'/y=u$ admits an algebraic solution~$y$.
\end{enumerate}
Step (ii) is solved by Risch's algorithm~\cite{Risch70}
(and independently by Baldassarri and Dwork~\cite[\S6]{BaDw79}); it was the missing ingredient in the 19th century procedures.
Singer solves Step~(i) by reducing it to an algebraic elimination problem. 
Part of the proof of the algorithm is that if step (ii) finds an
algebraic $y$ such that $y'/y=u$ with $u$ as in step (i), then $y$ is a nontrivial algebraic solution of $L$, and since $L$ is irreducible, all solutions of $L$ are algebraic.

In the reducible case, Singer's algorithm first proceeds with a factorization step that writes~$L$ as a product of irreducible operators. 
The corresponding detailed results and procedures can be found in (the proofs of)~\cite[Theorem~1]{Singer79} in the irreducible case, and~\cite[Theorem~3]{Singer79} in the general case. A slightly different route is proposed by Singer in~\cite[Corollary~4.3]{Singer81}, where first a basis of all Liouvillian solutions is computed, and then an algorithm by Rothstein and Caviness~\cite[\S4]{RoCa79} is used to decide if all elements in this basis are algebraic.

\subsection{Other approaches to problem {\bf (F)}}    \label{ssec:otherF}
Another way of seeing that Fuchs' problem~{\bf (F)} is decidable, in a
spirit similar to Singer's algorithm but of even higher algorithmic
complexity, 
is based on the equivalence for a linear differential operator~$L$ between having a full basis of algebraic solutions and having a finite differential Galois group.
Indeed, since Hrushovski's seminal article~\cite{Hrushovski02}, it is
known that one can fully determine algorithmically the Galois group of~$L$. (Prior to ~\cite{Hrushovski02}, it was only known how to compute the Galois group of \emph{completely reducible operators}, that is for operators $L$ that can be written as LCLM of irreducible operators~\cite{CoSi99}, a condition that can itself be algorithmically tested~\cite{Singer96}.)
The complexity of Hrushovski's algorithm is not yet fully understood,
and its simplification is the object of ongoing works, e.g., \cite{Feng15,Sun19,AmMiPo22}.
An approach to understanding the Galois group of $L$ by computing
the Lie algebra of its identity component is given in~\cite{BCDW20,BCDW16}
and~\cite{DreyfusWeil2022}.
In a different direction, van der Hoeven proposed 
in~\cite{Hoeven07} a symbolic-numeric approach for computing
differential Galois groups,
but for the moment its potential has not been exploited further.

To our knowledge, none of the algorithms mentioned in 
\cref{ssec:SingerF,ssec:otherF} has been implemented yet, nor are they expected to provide good practical behavior (except perhaps for very moderate orders, as in e.g., \cite{SiUl93b,Hessinger01}).

Let us finally mention a very recent sufficient (algorithmic)
criterion for problem~${\bf (F)}$.
Assume that the linear differential equation $L(y)=0$ admits a basis
of Puiseux series solutions at each of its singularities. Assume
moreover
that there exists another differential operator~$P$ such that
$\partial_z P$ is not right-divisible by $L$ and such that for all
singularities~$\sigma$ of~$L$ and for all elements~$f$ in a basis of
Puiseux
solutions of~$L$ at~$\sigma$, the Puiseux series $P(f)$ has only
nonnegative exponents.
Then $L$ admits a transcendental solution, in other words Fuchs' problem ${\bf (F)}$ is answered by the negative. This is the content of~\cite[Thm.~5]{KaKoVe23}, where such an operator $P$ is called a \emph{pseudoconstant}. Moreover, Algorithm~10 in~\cite{KaKoVe23} is able to compute a pseudoconstant for $L$, or to certify that none exists. 
More generally, if any symmetric power of $L$ admits a pseudoconstant,
then  $L$ admits a transcendental solution (Theorem~18 in~%
\cite{KaKoVe23}). Conversely, it is an open question whether the fact
that $L$ admits transcendental solutions implies the existence of
pseudoconstants for some symmetric power of~$L$ \cite[Question~20]{KaKoVe23}. If the answer to this open question
were positive,
then this would provide another way to show that Fuchs' problem is
decidable.

\section{Solving problem~{\bf (S)}} \label{ssec:newS} 

For a D-finite power
series $f \in \Q[[z]]$, we write $\Lmin$ for
the monic
linear differential operator in $\Q(z)\langle\partial_z\rangle$ 
of minimal order
that annihilates 
$f$. By a slight abuse of notation, we sometimes write~$L^{
\text{min}}_f$ for the operator in $\Q
[z]\langle\partial_z\rangle$ of the same order, obtained by clearing
denominators. These two operators differ by a multiplicative factor
in~$\Q(z)$. Both operators have the same main properties 
(solution space, local bases, singularities, \dots).
The  
\emph{degree} of $\Lmin$ is
the
maximal degree of its coefficients in the second representation.
$\Lmin$  can be computed efficiently~\cite{BoRiSa24}
starting from any linear differential equation satisfied by $f$
 together with sufficiently many initial terms of~$f$.

Our solution to Stanley's problem~{\bf (S)} is based on the properties of $\Lmin$ that we review now.

\subsection{Properties of \texorpdfstring{$\Lmin$}{Lmin
(f)}}\label{sub:properties_of_Lmin}

In simple terms, the minimal differential operator $\Lmin$ is a differential analogue for D-finite functions of the more classical notion of minimal polynomial for algebraic power series. 
Indeed, any $L\in\Q[z]\langle\partial_z\rangle$ such that $L(f(z))=0$
is a left multiple of $\Lmin$ in the ring $\Q
(z)\langle\partial_z\rangle$, and conversely any left multiple of 
$\Lmin$ is an annihilator for~$f$.

But there are two main differences between these similar notions. 
One is that, for such an $L$, although its order upper bounds the order of  $\Lmin$, it is in general not the case that the maximal degree of the coefficients of $L$ is an upper bound on the  maximal degree of the coefficients of
$\Lmin$; this is actually the main difficulty in the algorithmic computation of $\Lmin$, see~\cite{BoRiSa24} for more details.
Another difference important to note is that
$\Lmin$ 
does not need to be irreducible in $\Q(z)\langle\partial_z\rangle$.
Reducibility can already occur when $f$ is transcendental,
as illustrated by
$L^{\mathrm{min}}_{\log(1-z)} = ((1-z)\partial_z -1) \partial_z$. 
The same also holds when $f$ is algebraic.
For instance, for the algebraic power series of degree 2
\[
f=\sqrt{1-4 z}+z = 
1-z -2 z^{2}-4 z^{3}-10 z^{4}- \cdots, \]
we have that 
\begin{multline*}
\Lmin = 
\left(1 - 2 z \right) \left(1 - 4 z \right) \partial_z^{2}-
4 z \partial_z +4\\
=
\left(
\left(1 - 2 z \right) \left(1 - 4 z \right)  \partial_z +
4 z-6 + \frac{1}{z} \right)
 \cdot \left(\partial_z - \frac{1}{z} \right).
\end{multline*}
Moreover, this is actually the ``generic behavior'' for an algebraic power series $f\in\Q[[z]]$: if $P=y^d + c_{d-1}(z) y^{d-1} + \cdots + c_0(z) = \prod_{\ell=1}^d (y - f_\ell (z))$ is the minimal polynomial of $f=f_1$ in $\Q(z)[y]$, then $\Lmin$ annihilates all the conjugate roots $f_2, \ldots, f_d$ of $f=f_1$ (see~\cref{prop:basis_alg}). Therefore, $\Lmin$ admits $\sum_{\ell=1}^d f_\ell(z) = -c_{d-1}(z)$ as rational solution, hence if $c_{d-1} \neq 0$ then
$\Lmin$ is right-divisible by $\partial_z - c_{d-1}'/c_
{d-1}$.
This was remarked by Tannery \cite[p.~132]{Tannery1875},
see also \cite[Prop.~4.2]{CoSiTrUl02}.

\bigskip

In \cref{prop:basis_alg} below we give properties of $\Lmin$ and related operators.
The proofs rely on some basic facts from the Galois theory of linear
differential equations (see \cite{Singer90,PuSi03}), so we
begin in this
general
setting. Let $k$ be a differential field of characteristic zero with
derivation $\partial$ and whose field of constants  $C = \{ c \in k \
| \ \partial c = 0\}$ is algebraically
closed. Let $L \in k\langle
\partial\rangle$ be an operator of order $r$  
and let $K$ be the
associated Picard-Vessiot extension\footnote{This is a field
generated over $k$ by a fundamental set of solutions and all their
derivatives through order~$r-1$ and  having the same constants  $C$ as
$k$; one
sees that it is again a differential field. Since we assume that $C$
is
algebraically closed, such a  field exists and any two are
differentially isomorphic over $k$.}. Let $G$ be the differential
Galois group of $K$ over $k$ and let $f \in K$ be a solution of $L(y)
= 0$.  We define $\Lmin$ to be the monic operator in $k\langle
\partial\rangle$ of smallest order vanishing on $f$\footnote{This
definition is reconciled with the previous one before 
\cref{prop:basis_alg}.\hypertarget{fn:source}{}\label{fn:source}}.

The $C$-space of algebraic solutions of $L(y) = 0$ is left invariant
by the action of $G$, so it is the solution space of a monic operator
in $k\langle \partial\rangle$ denoted by $\Lalg$ (\cite[Lemma~2.2]
{Singer96} or \cite[Lemma~2.17 p.~48]{PuSi03})\hyperlink{fn:source}{\footnotemark[\getrefnumber{fn:source}]}.
Let $P(Y)$ be {an irreducible}
polynomial in $k[Y]$ and assume that there is a $u \in K$ such that $P
(u) = 0$. The differential Galois theory implies  that the splitting
field of $P$ over $k$ lies in $K$ \cite[Prop.~1.34.3]{PuSi03}.
Since $G$ permutes the roots of $P(Y)$, we have that the $C$-span of these roots is left invariant under this action. As before, this implies that this vector space is the solution space of a monic  operator in $k\langle \partial \rangle$, denoted by $\LP$. {Note that since Picard-Vessiot extensions for $L$ are unique up to $k$-differential isomorphisms, the operators $\Lalg, \Lmin$ and $\LP$ will be independent of such an extension since any such isomorphism will preserve all the properties of these operators.}

{We shall apply the above constructions to a linear differential
operator $L \in \Qbar(z)\langle \partial_z\rangle$. To be precise,
by an algebraic
solution we mean an element $y \in \calP = \bigcup_{n\in \NX^\star}
\Qbar
((z^{\frac{1}{n}}))$, the differential field of formal Puiseux series
where $\partial_z(z^{\frac{1}{n}}) = \frac{1}{n} z^{\frac{1}{n}-1}$
such that $y$ is algebraic over $\Qbar(z)$ and $L(y) = 0$. For a given
operator $L$, the set $\AL$ of algebraic solutions forms a vector
space over $\Qbar$ of dimension at most~$r=\ord(L)$.
Furthermore,
the  derivation $\partial_z$ extends uniquely to the field $E =
\Qbar(z, \AL)$. We now construct a Picard-Vessiot extension $K$ of
$k=\Qbar(z)$ that contains $\AL$. Let $F$ be the Picard-Vessiot
extension of $E$ corresponding to $L$. Since $\AL$ lies in the solution space of $L$ in $E$, there exists a basis $\{y_1, \ldots , y_r\}$ of this space that contains a basis of $\AL$. The field $K = \Qbar(z)(y_1, \ldots , y_r, \ldots, y_1^{(r-1)}, \ldots , y_r^{(r-1)}) \subset F$ has no new constants and so is the required Picard-Vessiot extension.}

If both~$k=\Qbar(z)$ and~$f\in\QX[[z]]$, 
\cref{prop:basis_alg} below shows that~$\Lmin$ and~$\Lalg$ have
coefficients
in~$\mathbb Q(z)$. Therefore in this special case, the definitions of
$\Lmin$ and $\Lalg$ coincide with those given earlier.

 \begin{proposition}\label[proposition]{prop:basis_alg}
  Let $ L \in \QX(z)\langle \partial \rangle$ and $f \in \QX[[z]]$ be
  an algebraic solution of $L(y) = 0$ with minimal polynomial $P(Y) \in \QX(z)[Y]$. If $\Lalg, \LP$, and $\Lmin$  are the operators defined above over $\Qbar(z)$, then these operators have coefficients in $\QX(z)$. Furthermore, $\Lmin = \LP$ so the solution space of $\Lmin$ is spanned by the roots of $P(Y)$. In particular, the local solutions of $\Lmin$ at any of its singularities do not contain logarithms.
\end{proposition}

\begin{proof}Once again let $\calP$  be the differential field of
Puiseux series and let $\Gal(\Qbar/\QX)$ be the Galois group of
$\Qbar$ over $\QX$. For $\sigma \in \Gal(\Qbar/\QX)$, we define the
action of $\sigma$ on $\calP$ as $\sigma(\sum_i a_i z^{\frac{i}{n}}) = \sum _i \sigma(a_i) z^{\frac{i}{n}}$. In this way $\Gal(\Qbar/\QX)$ acts on $\calP$ as a group of differential automorphisms with fixed field $ \bigcup_{n\in \NX^\star} \QX((z^{\frac{1}{n}}))$.

\smallskip
{$\boldsymbol\Lalg$}: Since $L$ has coefficients in $\QX(z)$,  $\Gal
(\Qbar/\QX)$ leaves the $\Qbar$-space of algebraic solutions
invariant. Denoting by $\sigma\Lalg$ the operator obtained from
$\Lalg$ by applying $\sigma$ to the coefficients of $\Lalg$, we have
that $\sigma\Lalg$ and $\Lalg$ have the same solution spaces and so must be equal. Therefore $\Lalg$ has coefficients in $\QX(z)$.

\smallskip
{$\boldsymbol\LP$}: Since $P(Y)$ has coefficients in $\QX(z)$, $\Gal
(\Qbar/\QX)$ permutes its roots and therefore  the $\Qbar$-space spanned by these roots is left invariant by this action. Arguing as above we have  that $\LP$ has coefficients in $\QX(z)$.

\smallskip
{$\boldsymbol\Lmin$}: Any $\sigma \in \Gal(\Qbar/\QX)$ leaves $f$
fixed
and sends $\Lmin$ to  $\sigma\Lmin$. By uniqueness of $\Lmin$, we have  $\Lmin = \sigma \Lmin$ and this implies, as above, that the coefficients of $\Lmin$ lie in $\QX(z)$.

\smallskip
{$\boldsymbol\Lmin \boldsymbol= \boldsymbol\LP$}: Let $F$ be the
splitting field of $P (Y)$ over
$\QX (z)$ and let $H$ be the (usual) Galois group of this extension. The derivation on $\QX(z)$ extends uniquely to a derivation on $F$ and the elements of $H$ commute with this derivation. Since $P(Y)$ is irreducible,  $H$ acts transitively on the roots of $P(Y)$. Since $\Lmin$ has coefficients in $\QX(z)$ and vanishes on one root of $P(Y)$, $\Lmin$ vanishes on all the roots and therefore on  the $\Qbar$-space spanned by the roots of $P(Y)$. This implies that $\LP$ is a right factor of $\Lmin$.  Minimality of $\Lmin$ and the fact that $\LP(f) = 0$ implies that $\Lmin$ is a right factor of $\LP$. Therefore these two operators are equal.

Finally, the local branches of the roots of $P$ are Puiseux series at
every finite point and at infinity. Since these roots span the
solution space of $\Lmin$, it admits at every point a
local fundamental system of Puiseux solutions and therefore no
logarithmic local solutions.
\end{proof}

\Cref{prop:basis_alg} implies  the following classical result (for
definitions of the terms see~\cite[Chap.~15-16]{Ince1956}):

\begin{corollary} \label[corollary]{prop:Fuchs}
The differential operator $\Lmin$ 
of an algebraic power series $f\in \Q[[z]]$
enjoys the following properties:
\begin{enumerate}
	\item it is Fuchsian;
	\item it  admits only rational exponents at all its singular
	points;
	\item the monic indicial polynomials at all of its singular points
	split in $\Q[z]$ into distinct linear factors.
\end{enumerate}
\end{corollary}
\begin{proof}By \cref{prop:basis_alg}, all solutions 
of~$\Lmin$ are algebraic. Fuchs proved that in this
situation,
the differential operator satisfies properties (1) and (2)~\cite
[\S6.II]{Fuchs1866}. Conclusion~(3) is a consequence of the fact that a multiple root of
 the indicial equation at $z=a$ always introduces logarithmic terms in
 the basis of local solutions at $z=a$~\cite[\S6.I]{Fuchs1866},
 see also~\cite[p. 405]
 {Ince1956}.
\end{proof}

\subsection{Examples}

\begin{exa}\label[example]{ex:Apery}
Here is a proof of the transcendence of Apéry's power series
\[
f(z) = \sum_{n \geq 0} A_n z^n, \quad \text{where} \; A_n = \sum_{k=0}^n \binom{n}{k}^2 \binom{n+k}{k}^2.
\]
 First, creative telescoping~\cite{Zeilberger91} produces a linear
 recurrence with polynomial coefficients satisfied by the sequence~$
 (A_n)_n$: $A_0 = 1$, $A_1 = 5$,
 \[(n+1)^3 A_{n+1} + n^3 A_{n-1} = \left( 2\,n+1 \right)  ( 17\,{n}^
 {2}+17\,n+5 ) A_n, \quad n\ge1.\]  
This recurrence is then converted into a linear differential equation 
$L(f) = 0$
satisfied by~$f(z)$, where $L$ is the differential operator
\[ 
L = (z^4-34z^3+z^2)\partial_z^3+(6z^3-153z^2+3z)\partial_z^2+(7z^2-112z+1)\partial_z +z-5 .
 \]  
 In the third step one certifies that  
 $\Lmin = L$, meaning that $L$ is already the minimal-order differential equation satisfied by $f(z)$. This can be done using the minimization procedure in~\cite[Algorithm 1]{BoRiSa24}\footnote{ 
 Note that in this very particular case, one could alternatively 
 prove that $L$ is actually irreducible (for instance, by showing that neither $L$, nor its adjoint, admit any non-trivial hyperexponential solutions).}.
 
 Finally, the indicial polynomial of $\Lmin$ at $z=0$ is $z^3$, and the transcendence of $f(z)$ follows from part (3) of \cref{prop:Fuchs}.
\end{exa}

The approach used in~\cref{ex:Apery} is the basis
of~\cref{algo:Stanley2a,algo:Stanley2b} below. 

Note that there are alternative proofs  
for the transcendence of the Apéry series~$f(z)$.
One of them relies on the combination of three (non-trivial) ingredients:
(i) the asymptotics
$A_n \sim \frac{(1+\sqrt{2})^{4n+2}}{2^{9/4} \pi^{3/2} n^{3/2}}$~\cite{McIntosh1996};
(ii) Flajolet's criterion (\cref{ssec:Flajolet}); 
(iii) the fact that $\pi$ is transcendental.
Our proof based on differential operators appears to be more ``natural'', since
proving the transcendence of a function should be easier than proving the transcendence of a number.
Another alternative proof uses~\cite[Thm.~5.1]{ChZu10} (see also~\cite
[\S6]{Zudilin08}); it relies on the fact that $f(z)$ admits modular
parametrizations, implying that, up to an algebraic transformation,
$f (z)$ coincides with the hypergeometric series
${_3F_{2}}\biggl[\begin{matrix}
  {1/2}\kern.707em {1/2}\kern.707em{1/2}\\ 
  {1}\kern.707em{1}
\end{matrix}\,;\, z\biggr]$, which is transcendental~(
\cref{ssec:hyp}).
This proof does not extend to more general Apéry-like series, such as
those
in Appendix~\ref{appendix:Apery}.

\begin{algorithm}[t]
\caption[]{Deciding transcendence of D-finite functions}
\label{algo:Stanley1}
\begin{algorithmic}[1]
\Require $L=a_r(z)\partial_z^r+\dots+a_0(z)$ with $a_i(z) \in \Q[z]$;
\Statex \qquad $f_0$: a truncated power series specifying
\Statex\qquad\qquad a unique solution~$f\in\mathbb Q[[z]]$
of $L (f)=0$.
\Ensure{Either \textsf{T} if $f$ is transcendental, or \textsf{A} if $f$ is algebraic.}
\State $\Lmin \coloneq$\textsf{MinimalRightFactor}{($L$,
$f_0$)}\Comment{Bostan-Rivoal-Salvy~\cite[Algorithm~1]{BoRiSa24}}
\If{$\Lmin$ has a basis of algebraic solutions}\Comment{Singer's algorithm~\cite{Singer79}}
\State{ ${\textsf{B}\coloneq\textsf{A}}$} 
\Else{ $ {\textsf{B}\coloneq\textsf{T}}$}
\EndIf\\
\Return $\textsf{B}$
\end{algorithmic}
\end{algorithm}   

\medskip
The next example is of a more combinatorial flavor. It exhibits a case
of a nontrivial D-finite power series $f(z)$ in $\Q[[z]]$ whose minimal-order operator $\Lmin$ is not only reducible, but (almost) split.

\begin{exa}\label[example]{ex:trident}
The generating function of the ``trident walks
in the quarter plane'' is the power series
	\[ f(z) = \sum_{n\ge0} a_n z^n = 1+2 \, z+7 \, z^2+23 \, z^3+84 \,
	z^4+301 \, z^5+1127 \, z^6+\cdots.\]          
Its $n$th coefficient $a_n$ counts all the
$	
	\raisebox{1.5mm}{\begin{tikzpicture}[scale=.3, baseline=(current bounding box.center)]
		\foreach \x in {-1,0,1}
		  \foreach \y in {-1,0,1}
	    \fill(\x,\y) circle[radius=1.5pt];
	{\draw[thick,->](0,0)--(1,1);}
	{\draw[thick,->](0,0)--(0,1);}
	{\draw[thick,->](0,0)--(-1,1);}
	{\draw[thick,->](0,0)--(0,-1);}
	\end{tikzpicture}}
$-walks of length $n$ in $\N^2$ starting at $(0,0)$.   
It was proved in~\cite{BoChHoKaPe17} that $f(z)$ satisfies $L(f) = 0$, where $L$ is a linear differential operator of order~5 and degree at most 15. 
In~\cite[\S4.2]{BoChHoKaPe17} it was shown that $f(z)$ is transcendental by exploiting the full factorization of~$L$. 
An alternative proof of the transcendence of~$f$ (and actually of all
$4\times 19 - 4 = 72$ similar power series in~\cite[Thm.~2]{BoChHoKaPe17}, see~\cref{sec:implem})
uses minimization to first prove that the differential equation $L
(f)=0$ has minimal order, that is $\Lmin = L$, and
concludes using the presence of logarithmic terms in the local
fundamental system of $L$ at $z=0$.
One advantage of this proof is that factorization of linear differential operators is avoided.
\end{exa}

Flajolet's asymptotic criterion used in the alternative proof for the
Apéry series does not prove transcendence in~\cref{ex:trident}.
Indeed, the asymptotic behavior\footnote{This
asymptotic estimate was conjectured in~\cite[Conjecture~1]{BoChHoKaPe17},
shown to be equivalent to an integral evaluation in~\cite[Theorem~8]{BoChHoKaPe17},
and proved using analytic combinatorics in several variables in~\cite[\S6]{MeWi19}.} 
$a_n \sim \gamma \beta^n n^r$ 
with
$\gamma = 4/(3\sqrt\pi)$, $\beta=4$, $r=-1/2$,
is compatible with algebraicity,
since 
$r \in \Q\setminus\Z_{<0}$, $\beta \in \Qbar$ and
$\gamma \Gamma(r+1) = 4/3\in \Qbar$. 

Another combinatorial example, also arising from the world of lattice path enumeration,
was considered in~\cite[Proposition~8.4]{BoBoKaMe16}.
In that case, all known criteria for transcendence fail to apply, as
well as all previously known algorithms.
Moreover, the defining differential operator (of order 11 and
polynomial degree 73) could not be factored with the computer
algebra systems of the time.
The proof in~\cite[Proposition~8.4]{BoBoKaMe16} relies on the
computation of the corresponding minimal operator~$L_f^{\text{min}}$.
This example was the starting point for 
\cref{algo:Stanley1,algo:Stanley2a,algo:Stanley2b}.

\subsection{An algorithm for Stanley's problem~{\bf (S)}} 
Based on \cref{prop:basis_alg}, a first solution to Stanley's problem~%
{\bf (S)} is \cref{algo:Stanley1}. It solves it \emph{in principle}, in the sense that it
shows that {\bf (S)} is decidable. However, its computational
complexity is too high to be of practical use, mainly because it
relies on the costly algorithm from~\cite{Singer79}. 

\begin{theorem} 
	Given a power series~$f\in\Q[[z]]$ 
		by a linear differential operator~$ L\in\Q[z]\langle\partial_z\rangle$ such that~$L(f)=0$, together with sufficiently many initial conditions,
\cref{algo:Stanley1} decides whether $f$ is algebraic or transcendental. 
\end{theorem}   

\begin{proof}
If $f$ is algebraic then, by \cref{prop:basis_alg},
$\Lmin$ admits a full basis of algebraic solutions, and
this is detected by Singer's algorithm~\cite{Singer79}. Conversely, 
if $\Lmin$ admits a full basis of algebraic solutions
then $f$ is a linear combination of the elements of this basis and
therefore is algebraic.
\end{proof}

\subsection{An efficient transcendence test} To go further towards a more efficient variant of~\cref{algo:Stanley1}, we use \cref{prop:basis_alg,prop:Fuchs}.

\begin{algorithm}[t]
\caption[]{Transcendence test for D-finite functions}\label{algo:Stanley2a}
\begin{algorithmic}[1]
\Require $L=a_r(z)\partial_z^r+\dots+a_0(z)$ with $a_i(z) \in \Q(z)$;
\Statex\qquad $f_0$: a truncated power series specifying
\Statex\qquad\qquad a unique solution~$f\in\mathbb Q[[z]]$ of $L(f)=0$.
\Ensure{Either \textsf{T} if $f$ is transcendental, or 
\textsf{FAIL}.}
\Statex 
\State $\Lmin\coloneq ${MinimalRightFactor}{($L$,
$f_0$)}\Comment{Bostan-Rivoal-Salvy algorithm~\cite{BoRiSa24}}
\For{$s\in\operatorname{Singularities}(\Lmin)$}
\State $P\coloneq \operatorname{IndicialPolynomial}(\Lmin,z=s)$
\If{$\deg P<\ord (\Lmin)$}{\ \Return $\textsf
T$}\label{line:line4}\Comment
{not Fuchsian}
\EndIf
\State $Z\coloneq \operatorname{DistinctRationalZeros}(P)$
\If{$\operatorname{card} Z<\deg P$}{\ \Return $\textsf
T$}\label{line:line6}\Comment
{the indicial polynomial does not have distinct rational roots}
\EndIf
\State $D\coloneq \{|z_i-z_j|,z_i\in Z,z_j\in Z\}\cap \mathbb N_{>0}$
\If{$D\neq\emptyset$}\label{algo:line8}
\State $B\coloneq \operatorname{FormalSolutions}(\Lmin,z=s,\mathrm{order}=\max(D))$
\If{$B$ has logarithms}{\ \Return $\textsf T$}
\EndIf
\EndIf
\EndFor
\State\Return \textsf{FAIL}
\end{algorithmic}
\end{algorithm}

\begin{theorem} 
	Given a power series~$f\in\Q[[z]]$ 
		by a linear differential operator~$ L\in\Q[z]\langle\partial_z\rangle$ such that~$L(f)=0$, together with sufficiently many initial conditions,
\cref{algo:Stanley2a} either proves that
$f$ is transcendental or returns {\sc FAIL}.
\end{theorem}   
\begin{proof}
Assume that $f$ is algebraic. By \cref{prop:Fuchs},
$\Lmin$ is Fuchsian and, at each of its singularities,
its indicial polynomial has degree $\ord (\Lmin)$ and has
distinct rational roots. Therefore neither of the tests in
lines~\ref{line:line4} and~\ref{line:line6} can return $\mathsf T$.

Moreover, by \cref{prop:basis_alg}, the local solutions of
$\Lmin$ contain no logarithmic terms. Thus the
logarithm test cannot return $\mathsf T$ either. Consequently,
\cref{algo:Stanley2a} returns $\mathsf{FAIL}$ whenever $f$
is algebraic. Hence, whenever it returns $\mathsf T$, the power
series $f$ is transcendental.
\end{proof}

In terms of completeness, \cref{algo:Stanley2a} is not as good as 
\cref{algo:Stanley1}, since it does not detect algebraicity. As such, it is not a full decision procedure for transcendence.
However, in terms of efficiency, \cref{algo:Stanley2a} is much better than 
\cref{algo:Stanley1}, as it only relies on the minimization algorithm from \cite{BoRiSa24}, plus a few computationally cheap local tests.

\subsection{An efficient conditional algorithm in the globally bounded case} 

To state our next results, we now recall a collection of
results and conjectures on D-finite globally bounded power series 
(see \cref{def:gb}). An important source of globally bounded power
series is the class of \emph{diagonals of rational functions} defined as
follows.
\begin{definition}For a rational function~$R
\in\mathbb Q(x_1,\dots,x_s)$ with Taylor expansion
\[\sum_{i_1\ge0,\dots,i_s\ge0}{a_{i_1,\dots,i_s}x_1^{i_1}\dotsm x_s^
{i_s}}\]
in the neighborhood of the origin, the diagonal~$\Diag(R)$ is the
power series
\[\Diag(R)=\sum_{i\ge0}{a_{i,\dots,i}z^i}\in\mathbb Q[[z]].\]
A power series in~$\mathbb{Q}[[z]]$ is a \emph{diagonal of a rational
function} if it
equals~$\Diag(R)$ for some rational function~$R$.
\end{definition}
Algebraic functions are diagonals~\cite{Furstenberg67}.
Diagonals of rational functions are D-finite and globally
bounded~\cite[Proposition~5]{Christol90}. Moreover, a D-finite
and globally bounded
power series is a $G$-function. 
A consequence is an extension of parts~(1)
and~(2) of \cref{prop:Fuchs}: $\Lmin$ is Fuchsian with rational
exponents at all its singular points. This follows from results of
Katz (1970), Chudnovsky and
Chudnovsky (1985) and André (1989), see e.g., \cite[p.~719]{Andre00}
or~\cite[Thm.~1]{lepetit} and the references therein\footnote{Note that part~(3) of 
\cref{prop:Fuchs} is
false
for the
larger class of diagonals, as the simple example 
$f 
= \Diag\left( {1}/(1-(x_{1}+x_2)(x_{3}+1))\right)
= \sum_{n \geq 0}  \binom{2n}{n}^2 z^n
$ 
demonstrates, since 
$ \Lmin = \left(16 z^{2}-z \right) \partial_z^{2}+\left(32 z -1\right) \partial_z +4
$
has indicial polynomial $z^2$ at
$z=0$.}.

The converse of some of these implications is the object of the
following conjecture. Notably, part~(3) extends part~(2) and
gives a converse to \cref{prop:basis_alg} that allows (conjecturally)
to detect algebraic solutions of linear differential equations under
the global boundedness assumption.

\begin{conjecture}[Christol, André]\label[conjecture]{conj:ChristolAndre}
Let $f\in \Q[[z]]$ be a globally bounded and D-finite power series. Then:
\begin{enumerate}
    \item $f$ is the diagonal of a rational function;
    \item If $z=0$ is an ordinary point for $\Lmin$, then $f$ is algebraic;
    \item If $\Lmin$ admits a logarithm-free fundamental
    system at~$z=0$, then $f$ is algebraic.
\end{enumerate}
\end{conjecture}

Part~(1) of \cref{conj:ChristolAndre} is due to Gilles
Christol, and was formulated in the late 1980s~\cite{Christol86}, see
also~\cite[Conj.~4]{Christol90} and~\cite[Conj.~10]{Christol15}.
Part~(2) was formulated around 1997 in private discussions between Gilles Christol and Yves André, and appeared in print in~\cite[Rem.~5.3.2]{Andre04} (there, the connection with the Grothendieck-Katz $p$-curvature conjecture is also discussed).
See also~\cite[\S3.3]{Matzat06} for a similar conjecture (called ``Eisenstein's Algebraicity Criterion'' by the author).
Part~(3) is due to Yves André (private communication).

The global boundedness assumption in \cref{conj:ChristolAndre} is crucial for the three parts.
For instance, without this assumption, Part (1) does not hold for $f(z) = \log(1-z)$, Part~(2) does not hold for $f(z)=\exp(z)$ and Part (3) does not hold for the power series 
$f(z) =  \pFq21{1/6 \ 5/6}{7/6} z$ considered in~\cite[Prop.~2.3]{KuJi03}, which actually is the particular case $(a, b) = (1/6, 5/6)$
of the following example.

\begin{exa} \label{ex:algsing}
There exist D-finite transcendental power series 
whose minimal-order operator 
is Fuchsian and such that all singular points are algebraic (not
logarithmic) singularities.
Take $a$ and $b$ in $\Q\setminus\Z$ such that $b-a$ is not an integer.
Then $f(z) =  \pFq21{a \ b}{a+1} z$ is such that its $\Lmin$ is reducible, and equals
\[
 z(1-z)\partial_z^2 +
(a+1-(a  +  b  +  1)z)\partial_z- ab = 
\Big( z (1 - z) \partial_z + 1 - (b  + 1) z \Big) \Big(\partial_z+a/z \Big) . \]
Thus, $\Lmin$ is Fuchsian, it admits the algebraic solution $z^{-a}$,
and its local exponents are 
$\{0,-a \}$ at $z=0$,
$\{0,1-b \}$ at $z=1$,
$\{a, b \}$ at $z=\infty$.
Since the exponent differences are not integers, the three singularities are algebraic.
Moreover, $f(z)$ is transcendental, since otherwise the operator $\Lmin$ would have finite cyclic monodromy 
and this is impossible by \cite[Thm.~2.2~(2)]{Vidunas07}.
However, $f(z)$ is not globally bounded by~\cite[Prop.~1]{Christol86}, so this example does not contradict Part~(3) of \cref{conj:ChristolAndre}. 
\end{exa}

\medskip                     
We design now a more efficient version of \cref{algo:Stanley1};
this is
\cref{algo:Stanley2b}, based on
\cref{prop:basis_alg,prop:Fuchs,conj:ChristolAndre}.

\begin{theorem}
	Given a linear differential operator~$ L\in\Q
	[z]\langle\partial_z\rangle$ and sufficiently many initial
	conditions specifying a unique solution~$f\in\Q [
	[z]]$, every output~\emph{\textsf{T}} of
\cref{algo:Stanley2b} is unconditionally correct. If $f$ is globally
bounded, then, assuming Part~(3) of 
\cref{conj:ChristolAndre}, every output~\emph{\textsf{A}} is also
correct.
\end{theorem}
\begin{proof}
If $f$ is algebraic, then, by~\cref{prop:Fuchs},
$\Lmin$ is Fuchsian and it admits only rational exponents at all its singular points. 
Moreover, if 
$\Lmin$ has a logarithmic singularity at $z=0$, then $f$ is
transcendental by \cref{prop:basis_alg}. This proves the
first part of the result.
The last part follows from \cref{conj:ChristolAndre}.
\end{proof}

\begin{rk}
	\Cref{ex:algsing} shows that the assumption ``\emph{$f(z)$ is globally bounded}''
is necessary in \cref{algo:Stanley2b}, even if we replaced  ``$\Lmin$ has a logarithmic singularity at $z=0$'' by
``$\Lmin$ has a logarithmic singularity at one of its singular points''.
\end{rk}

\begin{rk}\label{rk:diagonal_grade}
As already mentioned in \cref{ssec:criteria}, it is currently not
known how to decide if a given D-finite power series $f(z)$ is
globally bounded. However, this property is known to be satisfied by
large
classes of power series, such as generating functions of 
multiple binomial sums or, equivalently, by diagonals
of rational functions, see~\cite[Theorem 3.5]{BoLaSa17}. In this
case, the operator $L$ annihilating $f$ does not need to be part of
the input: it can be computed from a diagonal representation of~$f$ 
using \emph{creative telescoping}, see~\cite{BoLaSa17} and the references therein (\cref{sec:examples-diagonals} provides several non-trivial examples).
Still in this case,  \cref{algo:Stanley2b} can be enhanced in the following way: after computing $\Lmin$, if it detects logarithms in the local basis of solutions of $\Lmin$ at $z=0$, then it concludes transcendence; if it detects a term  $\log^{s}(z)$ in that local basis (e.g., if the indicial polynomial of $\Lmin$ at~$z=0$ is~$z^{s+1} \cdot P(z)$ with $P(0) \neq 0$), then it can conclude that the  number of variables needed to express $f(z)$ as the diagonal of a rational function is at least~$s+2$~\cite[Corollary~2.7]{HarderKramerMiller2026},
thus proving a strong form of transcendence of $f(z)$ if $s\ge1$. 
For
instance, in \cref{ex:Apery}, the presence of terms $\log^2(z)$ in the local basis of $\Lmin$ at $z=0$ implies that this number is at least~$4$ (and in fact exactly~$4$, e.g., using the diagonal representations in~\cref{sec:generic_diags}), see also~\cite[Example 3.9]{HarderKramerMiller2026}.
\end{rk}

\begin{algorithm}[t]
\caption[]{Transcendence of globally bounded D-finite functions}\label{algo:Stanley2b}
\begin{algorithmic}[1]
\makeatletter
\setcounter{ALG@line}{7}
\makeatother
\Require $L=a_r(z)\partial_z^r+\dots+a_0(z)$ with $a_i(z) \in \Q
(z)$;
\Statex\qquad $f_0$: a truncated power series specifying
\Statex\qquad\qquad a unique solution~$f\in\mathbb Q[[z]]$ of $L(f)=0$.
\Statex \qquad\emph{It is assumed that $f$ is globally
bounded.}
\Ensure{Either \textsf{T} if $f$ is transcendental, or \textsf{A} if
$f$ is algebraic (under \cref{conj:ChristolAndre}).}
\Statex 
\Statex \emph{In \cref{algo:Stanley2a}, replace \textsf{FAIL} by 
\textsf{A}
and line~\ref{algo:line8} by}
\If{$s=0$ and $D\neq\emptyset$}\EndIf
\end{algorithmic}
\end{algorithm}

\section{Solving problem {\bf (P)}}    \label{ssec:SingerP}
Given an operator $L \in \mathbb{Q}(z)\langle \Dz \rangle$ the
aim of this section is to present \cref{algo:Lalg} that
computes the
operator $\Lalg\in \mathbb{Q}(z)\langle \Dz \rangle$
whose solution space is spanned by the algebraic solutions of $L(y)=0$ 
(see \cref{sub:properties_of_Lmin}). 
This is done by extending some of Singer's 
results~\cite{Singer79,Singer96}.
\Cref{algo:Lalg} is then used in
\cref{algo:Stanley3} to
solve Stanley's problem {\bf (S)} in full generality, providing a second proof that Stanley's problem {\bf (S)} is decidable.  

Both algorithms rely on Singer's algorithm~\cite{Singer79}, which
suffers from a
very
large bound, except for very small order, where one can use more
refined bounds~\cite{Kovacic86,Hessinger01,SiUl93}, see 
\cref{ssec:SingerF}. This prevents us
from obtaining practical versions of \cref{algo:Lalg,algo:Stanley3}.
Still, the computation of $\Lalg$ is a question of
independent interest. 
We first recall the following definition:
\begin{definition}\label{def-equiv} Let $k$ be a differential field
of characteristic zero
and $L_1, L_2 \in k\langle \partial \rangle$  be {irreducible}
operators of order~$n$.  We say  $L_1$ is \emph{equivalent} to $L_2$ if there exist $a_0, a_1, \ldots , a_{n-1} \in k$, not all zero,  such that $L_1$ divides 
$L_2 \circ (a_{n-1}\partial^{n-1} + \cdots + a_0)$ on the right.
\end{definition}
This is a special case of a more general definition\footnote{Of
operators 
\emph{of the same type}, whose definition does not require
irreducibility.} but we
will only
need this version and the fact that
this is an equivalence relation~\cite[Cor.~2.6]{Singer96}. If $K$ is a
Picard-Vessiot extension
of $k$ containing
a full set of solutions of $L_1(y)=0$ and of $L_2(y)=0$ then the map
$\phi:y \mapsto a_{n-1}\partial^{n-1}y + \dots + a_0y$ maps the solution
space of  $L_1(y)= 0$ to the solution space of $L_2(y)=0$. Since $L_1$ and
$L_2$ are irreducible, this map is a bijection~\cite[Cor.~2.6]
{Singer96}. All
solutions of
$L_1(y)=0$ are algebraic then so are their images by~$\phi$ and
therefore the same is true for~$L_2(y)=0$. The converse is true, since
 the inverse of the bijection is of the same type.

\begin{algorithm}[t]
\caption[]{Algebraic Right Factor}
\label{algo:Lalg}
\begin{algorithmic}[1]
\Require{$L\in \Qbar(z)\langle \partial \rangle$.}
\Ensure{$\Lalg\in \Qbar(z)\langle \partial \rangle$
with solution space spanned by the algebraic solutions of~$L(y)=0$.}
\Statex
\State\label{P-step1} factor $L$ over $\Qbar(z)$ as a
product $L_1 L_2
\cdots
	L_t$ where each $L_i$ is irreducible
\State\label{P-step2} select from the factors $L_i$ those factors $L_{i_1},
\ldots, L_{i_s}$\par that admit a full basis of algebraic
solutions\Comment{Singer's algorithm~\cite{Singer79}}
\State\label{P-step3}
for $j = 1, \ldots s$, construct an operator $\widehat{L}_{i_j}$
such that\Comment{See \cref{lem:2}}
\[ \widehat{L}_{i_j} = \textrm{LCLM}\{ R \in \Qbar(z)\langle \partial
\rangle \ | \ R \text{ is equivalent to } L_{i_j} \text{ and  divides } L \text{ on the right}\}\] 
\State compute $\Lalg=\operatorname{LCLM}\{\widehat{L}_
{i_1},\dots,\widehat{L}_{i_s}\}$\label{P-step4}
\State \Return $\Lalg$
\end{algorithmic}
\end{algorithm}        

\medskip

\begin{algorithm}[t]
\caption[]{Deciding transcendence of D-finite functions}
\label{algo:Stanley3}
\begin{algorithmic}[1]
\Require $L=a_r(z)\partial_z^r+\dots+a_0(z)$ with $a_i(z) \in \Q[z]$;
\Statex \qquad $f_0$: a truncated power series specifying
\Statex\qquad\qquad a unique solution~$f\in\mathbb Q[[z]]$
of $L (f)=0$.
\Ensure{Either \textsf{T} if $f$ is transcendental, or \textsf{A} if $f$ is algebraic.}
\State $\Lalg \coloneq$\textsf{AlgebraicRightFactor}{($L$)}\Comment{\cref{algo:Lalg}}
\If{$\Lalg(f_0)=0$}{ ${
\textsf{B}\coloneq \textsf{A}}$} \Comment{\cite[Lemma~2.1]{BoRiSa24}}
\Else{ $ {\textsf{B}\coloneq\textsf{T}}$}
\EndIf
\State\Return $\textsf{B}$
\end{algorithmic}
\end{algorithm}

\Cref{algo:Lalg} computes $\Lalg$ 
when  $k=\Qbar(z)$ and $C=\Qbar$.
Step~\ref{P-step1} is performed using any factorization algorithm, 
e.g., \cite{Grigoriev90}, or~\cite{Hoeij1997a}; step~\ref{P-step2}
is solved
by Singer's algorithms~\cite{Singer79} or \cite{Singer81}; for the
computation of the LCLM in
step~\ref{P-step4}, see e.g., \cite{BoChLiSa12}, or~\cite{Hoeven16}.
By \cref{prop:basis_alg}, if the input~$L$ is in $\mathbb Q
(z)\langle\partial\rangle$, then even if the intermediate computations
take place in $\Qbar(z)\langle\partial\rangle$, the output $\Lalg$ does belong
to~$\mathbb Q(z)\langle\partial\rangle$.

We will now justify the equation in the last step and
show how one can
construct the operator in step~\ref{P-step3}.

\subsection{Proof of the algorithm}\label{sec:proof-algo}
This proof uses again concepts from the differential Galois theory of
linear
differential equations as in \cref{sub:properties_of_Lmin}.
 Let $k$ be a differential field of characteristic~0 with
 algebraically closed field of constants~$C$, $L \in k\langle
 \partial \rangle$,
 $K$ the Picard-Vessiot extension of $k$ corresponding to $L(y)=0$ and
 $G$ the Galois group of $K$ over~$k$.

We recall that the vector space of solutions of $L(y)=0$ that are algebraic over~$k$ is
invariant under $G$ and so spans the solution space $V^{\textrm{alg}}$ of a linear differential equation $\Lalg (y)=0$ with $\Lalg\in k\langle \partial \rangle$.
 Note that $\Lalg$ divides $L$ on the right  \cite[Lemma
 2.1]{Singer96}. The equation $\Lalg (y) = 0$ obviously has a
 full set of solutions in $K$ and they generate a Picard-Vessiot extension $E$ over $k$ with $k \subset E \subset K$.  Furthermore the Galois group $H$ of $E$ over $k$ is finite.  A finite group acting on a vector space allows one to decompose the vector space as a direct sum of irreducible subspaces.  
Therefore $V^{\textrm{alg}} = \oplus_i V_i$, each $V_i$ an
$H$-irreducible space. Furthermore, each $V_i$ is the solution space
of an irreducible $M_i \in  k\langle \partial \rangle$. Factorizations
of linear operators are unique up to equivalence~\cite[Proposition
2.11]{Singer96}. It follows that $M_i$ must be equivalent to some $L_
{i_j}$ and so must divide some  $\widehat{L}_{i_j}$. Therefore the solution space of  $\textrm{LCLM}\{\widehat{L}_{i_1}, \ldots , \widehat{L}_{i_s}\}(y)=0$ contains the solution space of $\Lalg (y) = 0$. Conversely all solutions of $\textrm{LCLM}\{\widehat{L}_{i_1}, \ldots , \widehat{L}_{i_s}\}(y)=0$ are algebraic so its solution space is contained in the solution space of $\Lalg (y) = 0$. Therefore these monic operators are equal.

\subsection{Computing the \texorpdfstring{$\widehat{L}_{i_j}$}{Lij}}
This will  follow from 
\begin{proposition} \label{lem:2}
Let $k$ be a differential field of characteristic~0 with
 algebraically closed field of constants~$C$ and satisfying the
 following property:
\begin{quotation} \noindent If $L \in k\langle \partial \rangle$ then
one can effectively find a $C$-basis of the space \[\{u \in k \ | \ Lu =
0\}.\]  \end{quotation}
Let $M, L \in k\langle \partial \rangle$ and assume that $M$ is
irreducible. One can effectively find\footnote{This is a LCLM
of a potentially infinite set of operators that are all right
divisors of $L$. Thus, the earlier definition extends and the LCLM
is a right divisor of~$L$. By convention, if the set is empty, the
LCLM is taken to
be~1.}
\[\widehat{M} = \textrm{\emph{LCLM}}\{ R \in k\langle \partial \rangle \ | \ R \text{ is equivalent to } M \text{ and  divides } L \text{ on the right}\}.\]
\end{proposition}
Note that in the application of this proposition in \cref{algo:Lalg},
where the field~$k$ is~$\Qbar(z)$, the property is satisfied by
Liouville's algorithm~\cite{Liouville1833}, see also~\cite{AbKv91}.
\begin{proof} Let $M$ have order $n$. The set 
\[ W= \{\mathbf{a} = (a_{n-1}, \ldots , a_0)  \in k^n \ | \ M \text{ divides } L \circ (a_{n-1}\partial^{n-1} + \dots + a_0) \text{ on the right}\}\]
is a vector space over the constants.   
By computing the product~$L\circ (a_{n-1}\partial^{n-1} + \dots +
a_0)$, its right Euclidean division by~$M$ and writing that the
remainder is~0, we get an $n\times n$ matrix $\calA$ with entries in
$k\langle \partial \rangle$ such that $\mathbf{a} \in W$ if and only if $\calA \mathbf{a} = 0$ (see the discussion
following Example 2.8 in \cite{Singer96}).  Finding solutions of
$\calA Y = 0$ in~$k^n$ can be reduced to finding solutions of scalar
linear differential equations in $k$ and so one can find a $C$-basis 
$\{\ba_1, \ldots , \ba_s\}$
of~$W$~\cite[Chap.~III, \S11]{Poole1960}.
If $W=\{0\}$, there is no~$R$ equivalent to~$M$ that divides~$L$ on
the right and the result is~$\widehat{M}=1$.

For $\ba = (a_{n-1}, \ldots , a_0) \in W\setminus\{0\}$,
we construct an operator
$L_{\ba} \in k\langle \partial \rangle$ of order at most~$n$ so that
the map $A_\ba:y \mapsto a_{n-1}\partial^{n-1}y + \dots + a_0y$
maps the
solution space of  $M$ to the solution space of $L_\ba$.
This is done by differentiating $j$~times the equality $v=a_
{n-1}\partial^ {n-1}y
+
\dots + a_0y$ for $j=0,\dots,n$ and using the equation $My=0$ to
express
each derivative~$v^{(j)}$ in terms of $y, \ldots , \partial^{n-1}y$.
The
resulting $n+1$ expressions must be linearly dependent over $k$ and
this yields an equation~$L_\ba$. The operator~$L_\ba$ obtained that
way has order exactly~$n$ and is irreducible: by its construction,
there exists an operator~$Q_1$ such that~$L_\ba\circ A_\ba=Q_1\circ
M$ with $A_\ba$ relatively prime to~$M$ since~$M$ is
irreducible; the conclusion follows from general decomposition
theorems~\cite[Ch.~II. Thm.~1]{Ore1933}, \cite[Prop.~2.11]{Singer96}.
This shows that $L_\ba $ is equivalent to $M$. Moreover, $L_\ba$
divides $L$ on the right since $A_\ba $ maps the solution space of 
$M$ bijectively into the solution space of $L_\ba$, while
it also maps it into the solution space of $L$.

Let ${\widetilde{M}} = \textrm{LCLM}\{L_{\ba_1}, \ldots , L_
{\ba_s}\}$. We will show that $\widetilde{M} = \widehat{M}$. Since
each $L_{\ba_i}$ is equivalent to $M$ and divides $L$ on the right, we
have that ${\widetilde{M}}$ divides ${\widehat{M}}$ on the right.

Now assume that an operator $R$ is equivalent to $M$ and divides $L$
on the right. Let $A_\ba=a_{n-1}\partial^{n-1} + \dots + a_0$ be such
that
$M$ divides $R\circ A_\ba$ on the right. Since
$R$ divides $L$ on the right,  
we have that  $M$ divides
$L\circ A_\ba$ on the right. Therefore $\ba = 
(a_{n-1}, \ldots , a_0) \in W$ and so  there exist $c_1, \ldots ,c_s
\in C$ such that $\ba =  \sum c_i \ba_i$.
Moreover, $A_\ba$ maps the solution
space of~$M$ bijectively into the solution space of~$R$.
Since $A_\ba=\sum{c_iA_{\ba_i}}$, the latter space is included in the
sum
of the solution spaces of the $L_{\ba_i}$,
which is the  solution space of ${\widetilde{M}}$.
Therefore $R$ divides ${\widetilde{M}}$ on the right. 

Since ${\widehat{M}}$  is the $\textrm{LCLM}$ of all such $R$, we have
that ${\widehat{M}}$
divides ${\widetilde{M}}$ on the right.
Since both of these operators are monic, we have shown that they are
equal.
\end{proof}

\begin{rk}
 \Cref{lem:2} can be used as a building block to construct the Loewy
 decomposition of an operator into the product of completely reducible operators.
\end{rk}

\section{Other algorithmic approaches to algebraicity and
transcendence} \label{ssec:other}
\subsection{Using the Grothendieck-Katz \texorpdfstring{$p$}{p}-curvature conjecture} \label{ssec:GK}
In the 1970s, Alexander Grothendieck proposed a conjectural ``arithmetic'' characterization for linear differential operators with coefficients in $\Q(z)$: such an operator admits a full basis of algebraic solutions if and only if the same holds for its reductions modulo~$p$ for almost all primes~$p$.
Grothendieck's conjecture can be seen as a differential generalization
of a particular case, due to Kronecker, of Chebotarev's theorem: the
roots of a polynomial in $\Q[z]$  are all rational if and only if the roots of its reductions modulo $p$ are in $\mathbb{F}_p$ for almost all primes~$p$.
(See~\cite{BoCaRo24,FuPa26}, and the references therein.)

Although it remains largely open in the general case, the $p$-curvature conjecture has been proven in several important cases. It was first studied in depth and popularized by Nicholas Katz~\cite{Katz72}, who proved it for the Gauss hypergeometric equations, a result later extended to the generalized hypergeometric equations by Beukers and Heckman~\cite{BeHe89}.
More generally, it was proved in~\cite{Katz72} and~\cite[\S III]{Andre04} for Picard-Fuchs differential operators (satisfied by the periods of families of smooth algebraic varieties) and for certain of their direct factors and subquotients (sometimes called ``differential equations coming from geometry'').
In a different direction, the $p$-curvature conjecture was proved by the Chudnovsky brothers~\cite{ChCh85} (see also \cite[\S6]{Andre04}) for linear differential operators of first order with algebraic function coefficients in~$\overline{\Q(z)}$.
(This last arithmetic approach was used to prove the algebraicity of
some power series e.g., in~\cite{Kontsevich11,DeRi24,Bostan25}.)

Consider now problem~{\bf (F)} for any linear differential operator~$L$ for which the Grothendieck-Katz conjecture is known to hold.
For any fixed prime~$p$,
the linear differential operator $L \bmod p$, with
coefficients in $\mathbb{F}_p(z)$, admits a basis of algebraic
solutions if and only if it admits a basis of rational
solutions~\cite{Chambert-Loir2002}, and this is verified by a very simple algorithm,
based on checking the nullity of the $p$-curvature of $L$.

In practice, computing $p$-curvatures for a few primes $p$ provides a very efficient heuristic: if $L$ has several zero $p$-curvatures, then  
it very likely admits a full basis of algebraic solutions. 
In particular, for problem~{\bf (S)}, if the Grothendieck-Katz
conjecture is known to hold for $\Lmin$ and if $\Lmin$ has several
zero $p$-curvatures, then $f(z)$ is very likely algebraic, otherwise
it is very likely transcendental.
To turn this heuristic into a decision algorithm, one would need an
upper bound on the number of primes~$p$ for which the $p$-curvature is to be tested. In other words, one needs an \emph{effective} Grothendieck-Katz theorem, similar to effective versions of Chebotarev's theorem (e.g., \cite{PiTuWu20} and the references therein).
For Picard-Fuchs operators, the existence of such results
is alluded to by André~\cite[\S16.3.1]{Andre04}.
A first complete effective result was recently obtained, for first-order linear operators with rational function coefficients, by Fürnsinn and Pannier~\cite{FuPa26}. 
Much remains to be done in this direction, for more general classes of linear differential operators.

\subsection{Guess-and-prove} \label{ssec:GnP}
To prove algebraicity of a given power series $f(z)$, a very popular method in computer algebra is to first guess a polynomial $P(z,u) \in \Q[z,u]$ such that $P(z,f(z))=O(z^\sigma)$ for some value~$\sigma$, and then to certify that $P(z,f(z))=0$ by manipulating the roots of $P$ via the differential operator $L_P$ whose solution space is generated by these roots.
This strategy is explained for instance in the proof of~\cite[Corollary 2]{BoKa10}. The guessing part is based on (structured) linear algebra (over~$\Q$) in size~$\sigma$
and the most efficient algorithms use Hermite-Pad\'e approximants (which reduce the problem to linear algebra over~$\Q(z)$ in a size much smaller than~$\sigma$).
The proving part is based on the computation of~$L_P$, which
ultimately also relies on linear algebra~\cite{BoChLeSaSc07}, and
sufficiently many initial conditions.
A delicate issue is the choice of~$\sigma$; usually one keeps doubling
its value and degree bounds on~$P$ until a candidate~$P$ is found. 
If $f(z)$ is algebraic, then the procedure will eventually discover it (modulo computational difficulties related to the size of~$\sigma$ that can be huge). 

This approach applies to problem {\bf (S)} when~$f$ is D-finite.
Its strength is that it is able to recognize and prove
algebraicity of $f(z)$ even in situations where $f(z)$ is not \emph{a
priori} known to be D-finite (e.g., it is given not by an ODE, but by
a different kind of functional equation), see e.g., \cite{BoKa10}.
The weak point of the approach is that it cannot prove that $f(z)$ is transcendental. This is why it is usually used in conjunction with other heuristics such as computations of $p$-curvatures, or numerical computations of local monodromies.

\subsection{Factoring} \label{ssec:factor}
If the input $L$ is checked to be irreducible, then its solution $f$ is algebraic if and only if all solutions of $L$ are algebraic (by Proposition 2.5 in~\cite{Singer79}), and this can be tested using Singer's 1979 algorithm~\cite{Singer79}. 
But if $L$ is reducible, 
then it is \emph{a priori} not clear how to combine factoring algorithms~%
\cite{Hoeij1997a,Grigoriev90} and Singer's algorithm~\cite{Singer79} in order to solve Stanley's problem {\bf (S)}. 
One practical way is: factor the input $L = AB$ for irreducible~$A$ of order less than the order of~$L$ and check if $A$ has a basis of algebraic solutions.
If yes, then we cannot conclude. 
{(See the examples at the end of \cref{sec:examples-diagonals}.)}
If no, then $A$ has no nontrivial algebraic solutions, and there are two possible cases.
Either $g\coloneq B(f)$ is nonzero, in which case it is transcendental (as it is a nontrivial solution of $A$), hence $f$ itself is transcendental (otherwise $g$ would be algebraic).
Or, $g=0$, and we can repeat the whole procedure on~$B$.
This procedure is guaranteed to succeed in the particular case when $L$ factors as a product of irreducible operators each of them having no algebraic solutions.

\section{Examples} \label{sec:implem}

Starting from an implementation of the minimization algorithm~\cite{BoRiSa24},
the implementation of the
transcendence test of \cref{algo:Stanley2a}
in Maple is straightforward and fits in 40~lines.\footnote{It is
implemented in the function \texttt{istranscendental} in versions
later than 4.05 of the package \texttt{gfun}, available at
\url{https://perso.ens-lyon.fr/bruno.salvy/software/the-gfun-package/}, where a Maple session with these examples can also be found.} We list a few experiments performed with
this implementation.

\subsection{Combinatorial examples}
Examples with a combinatorial origin for which our
implementation gives an
automatic proof of transcendence are: 
\begin{table}
\begin{tabular}{*5r}
\makecell{dim}&\makecell{order}&\makecell{degree}&
\makecell{transc.\\ (sec.)}&\makecell{factor.\\ (sec.)}\\
\hline
3&3&5&2.2&0.9\\
4&4&10&0.1&0.2\\
5&6&17&0.6&1.3\\
6&8&43&4.2&18.3\\
7&11&68&24.0&265.\phantom{0}\\
8&14&126&174.9&4706.\phantom{0}\\
9&18&169&771.6&$>$10000\\
10&22&300&8817.1& --- \\
11&27&409&46743.& --- \\ \hline\\
\end{tabular}
\caption{Transcendence of lattice Green functions.}
\label{table:ex}
\end{table}
\begin{itemize}
\item the differential equation
 of order~11 and degree~73 from
3-dimensional walks confined to the positive octant mentioned after~
\cref{ex:trident}, coming from \cite[Proposition~8.4]{BoBoKaMe16} (in
10 seconds)\footnote{{Timings were obtained with
Maple2024.2 on a MacBook Pro 2017 with a quad-core Intel i7.}};
\item the 72~transcendental cases
from~\cite[Theorem~2]{BoChHoKaPe17} (in
a total of 25~seconds);
\item the lattice Green functions of the face-centered cubic lattice
in dimensions $d=2,\dots,11$. 
These (D-finite) functions are defined by 
\[G_d(z) \coloneq \frac{1}{\pi^d} \int_0^\pi \cdots \int_0^\pi \frac{d\theta_1 \cdots d\theta_d}{1-z\binom{d}{2}^{-1} \sum_{1\leq i < j \leq d} \cos \theta_i \cos \theta_j}
.
\]
When $d=2$ and $d=3$, we have hypergeometric expressions \cite[Eq.~(53)]{Guttmann10}, \cite[Eq.~(3.6)]{Joyce12}
\[G_2(z) = \pFq21{\frac12\ \frac12}1
{z^2},\quad G_3(z) = \pFq21{\frac16\ \frac13}1 {\frac{z^2(z+3)}
{4}}^2\]
that can be seen to be transcendental by the arguments of \cref{ssec:hyp}.
For $d \geq 4$, no such closed formulas are known, and it is not clear how to decide transcendence of $G_d(z)$ using asymptotic arguments.
Linear differential equations for $G_d(z)$ were first conjectured by Guttmann
in dimension~4~\cite{Guttmann2009a} and by Broadhurst in
dimension~5~\cite{Broadhurst2009}, after intensive computations. 
They were then computed by Koutschan for $d=3,4,5,6$ using creative
telescoping~\cite{Koutschan2013b}. 
For $d=7$, an equation was conjectured by Hassani, Maillard and Zenine~\cite{ZenineHassaniMaillard2015} and then equations were
conjectured up to $d=11$ by Hassani,
Koutschan, Maillard and Zenine~\cite{HassaniKoutschanMaillardZenine2016}.

While these differential equations
have been obtained after days of
computations, it is relatively simple for our code to prove that they are minimal. 
(Note that their minimality for $d\leq 11$, 
and even their irreducibility for $d\leq 7$, were claimed in~\cite
[Section~2.2]{HassaniKoutschanMaillardZenine2016} and~\cite[Section~4]
{ZenineHassaniMaillard2015}, but the
arguments
there are based on heuristics.)
Since all these equations are minimal and have a logarithmic singularity at~0,
\cref{prop:basis_alg} implies that $G_d(z)$ is transcendental for
$d\leq 11$ (assuming the
conjectured differential equations to be correct for $d\ge7$).
Alternatively, one can use factorization of linear differential
operators to prove that the differential equations for $G_d(z)$ are
irreducible for $d \leq 8$, and use Proposition~2.5 in~\cite{Singer79} (rather than \cref{prop:basis_alg}) to deduce the transcendence of $G_d(z)$.
Sizes and timings are reported in \cref{table:ex}. (The column
``transc'' indicates the time taken by our code to prove transcendence of $G_d(z)$;
the time taken by the factorization program 
from~\cite{ChGoMe22} is given in the column ``factor''. Maple's native
command \texttt{DFactor} was much slower.)
\end{itemize}

In practice, for all these examples, all the
time is spent in the first step where a differential
equation of certified minimality is computed.

\subsection{Diagonals}\label{sec:examples-diagonals}
A classical result by P\'olya and
Furstenberg states that the diagonal of any bivariate rational
function is algebraic~\cite{Polya22,Furstenberg67}. In
characteristic~0 and with three or more
variables, diagonals need not be algebraic.
The examples below illustrate how our methods can prove transcendence and can suggest algebraicity, which may then sometimes be established by a separate argument. 

\subsubsection*{Example~1: an algebraic case} Using
Koutschan's
\texttt{HolonomicFunctions} package~\cite{Koutschan2010},
the diagonal of 
\[\frac1{(1 - 5x-7yz - 13z^2)(1-x-xy)}\]
is found to satisfy a linear differential equation of order~5 with
coefficients of degree up to~24. This equation has logarithmic
singularities at~0 and~$1/140$. However, minimization shows that the
diagonal actually satisfies an equation of
order only~3 with coefficients of degree up to~13. From that equation,
\cref{algo:Stanley2b} strongly suggests that the diagonal is algebraic
(being a diagonal, it is globally bounded). Actually, the more general
situation of the diagonal of
\[F(x,y,z)=\frac1{(1 - ax-byz - cz^2)(1-x-xy)}\]
for arbitrary $a,b,c$ can be seen to be algebraic by a direct
residue computation: the diagonal is the residue at $x=y=0$ of 
\[G(x,y,z,t)=\frac1{xy}F(x,y,t/(xy))=\frac{xy}{(
(x-ax^2-bt)xy^2-ct^2)(1-x-xy)}\]
when $t\rightarrow0$. The second factor of the denominator
has a root in~$y$ that is away from~0 while the first
factor has two roots in~$y$ that both
tend to~0 when $t\rightarrow0$. The sum of their residues is therefore
a rational function in~$x,t$. From there it follows that the
sum of residues in~$x$ is algebraic. An explicit computation is
possible: by a Rothstein-Trager resultant~\cite[Thm. 2.1]{FuPa26}, one obtains a polynomial 
cancelling these two residues in~$y$, from whose coefficients their sum is
deduced to be
\[\frac{x -1}{a \,x^{4}-2 a \,x^{3}+b t \,x^{2}+c \,t^{2} x +a \,x^
{2}-2 b t x -x^{3}+b t +2 x^{2}-x}.\]
Only one of the roots of the denominator tends to~0 with~$t$. Again, a
Rothstein-Trager resultant gives a polynomial of degree~4 cancelling
the diagonal.

\subsubsection*{Example~2: a transcendental case}
A very similar-looking diagonal $\Delta(z)$ is the diagonal of
\[\frac1{(1 -x -y -z^2)(1-x-xy)}.\] 
\texttt{HolonomicFunctions} gives a linear differential operator~$L_7$
of order~7 for $\Delta(z)$ with coefficients of degree up to~19. This operator is not
irreducible: it is the product~$T_4A_3$ of an operator~$T_4$ that does
not have nonzero algebraic solutions\footnote{This can be proved by
showing that $T_4$ is irreducible and admits local logarithms at $z=0$.} 
with an operator~$A_3$ {all of
whose solutions are algebraic}\footnote{This can be checked by
observing that $A_3$ is irreducible and admits the solution 
$$\frac1z\cdot {}_3F_2
\left[\left.\begin{matrix}1/6,1/2,5/6\\
1/4,3/4\end{matrix}\right |\frac{27z^2}{256}\right],$$ which is algebraic
by the interlacing
criterion of \cref{ssec:hyp}.}.
Our code proves that~$L_7$ is minimal for
the diagonal~$\Delta$, which is then proved transcendental by showing
that~$L_7$ has a logarithmic singularity at~0. Another proof 
(as in \cref{ssec:factor}) is by
observing that~$A_3$ does not annihilate~$\Delta$, which implies that
$A_3(\Delta)$ is a solution of~$T_4$ and thus necessarily
transcendental, and therefore so is~$\Delta$. That alternative way does
not work on the adjoint~$L_7^*$ of $L_7$. If one takes the
solution~$s(z)$ {of~$L_7^*$} which is 
its unique power series solution with $s(0)=-1/2$, $s'(0)=0$,
$s''(0)=-21121726441/112000$, then again our program proves that it is
transcendental. In this example, the factorization~$L_7^* = A_3^*T_4^*$
has the left factor $A_3^*$ whose solutions {are all
algebraic}\footnote{A very similar hypergeometric function appears,
with
parameters $\{ 1/6, 1/2, 5/6 \}, \{ 5/4, 7/4 \}$.}, but now this
factorization is not sufficient to conclude that~$s(z)$ is
transcendental.

\subsubsection*{Example 3: a transcendental case with a difficult
right factor}

Finally, let us consider the diagonal of
\[\frac1{(1-x-y-z^2)(1-x-xy^2)} .\]
It is annihilated by a linear differential
operator~$L_9$ of order~9 and degree~60 that 
factors\footnote{This operator and its factorization were
communicated to us by Jean-Marie Maillard. Similar examples can be
found in~\cite{HassaniMaillardZenine2025}.}
as a product~$L_4A_5$. 
Our code detects that $L_9$ is minimal for the
diagonal and proves its transcendence. 
Another transcendence proof is based on the method in~\cref{ssec:factor}, as above.
By contrast, it is much more challenging to prove that all solutions
of~$A_5$ are algebraic.
The heuristic based on $p$-curvatures in~\cref{ssec:GK} suggests that this is indeed the case.
The approach using guess-and-prove described in~\cref{ssec:GnP} is applicable in principle: the guessing part is feasible in practice, although computationally very expensive, but the proving step seems to be out of reach, even modulo a prime\footnote{Gilles~Villard 
found a polynomial of degree~120
with coefficients of degree~460 that cancels the reduction
modulo~$p$ of a power series with
largest valuation among the solutions of~$A_5$, modulo a large prime
$p=892901$. Similar polynomials can be found for a basis of
solutions $\bmod p$.
Proving that the solutions of such a polynomial are a basis
of the reductions modulo~$p$ of the solutions of~$A_5$ is feasible in
theory, but is extremely demanding
computationally.}.

\subsection{Apéry-like series}\label{sec:Apery_like}
For $p,q \in \N\setminus \{ 0 \}$, we consider the following generalization of the power series in~\cref{ex:Apery}:
\[
f_{p,q} (z) =
\sum_{n \geq 0} \left( 
\sum_{k=0}^n 
\binom{n}{k}^{p} 
\binom{n+k}{k}^{q} 
\right) z^n.
\]
In \cref{appendix:Apery} we use Flajolet's criterion to prove that 
$f_{p,q} (z)$ is transcendental if and only if $(p,q) \neq (1,1)$.
It was conjectured by Sergey Yurkevich in~\cite[Conjecture~10.5]{Yurkevich2023} that the minimal number $\mu_{p,q}$ of variables needed to represent $f_{p,q}(z)$ as the diagonal of a rational power series is equal to $p+q$.
Yurkevich showed that $\mu_{p,q}\leq p + q$, by using the identity (due to Wadim Zudilin)
\begin{equation}\label{eq:Apery_as_diag}
f_{p,q} = \Diag \! \left( \frac{1}{\left(\prod_{j=1}^q (1-y_j)-x_1 \right) \cdot \prod_{k=2}^p (1-x_k) - \prod_{k=1}^p x_k \cdot \prod_{j=1}^q y_j} \right) .
\end{equation}
For $p+q \leq 10$, we checked the opposite inequality 
$p + q \leq \mu_{p,q}$ by using \cite[Corollary~2.7]{HarderKramerMiller2026} and the enhanced version of our \cref{algo:Stanley2b}, mentioned in~\cref{rk:diagonal_grade}. 
Indeed, in all these 45 cases, the indicial polynomial $I_{p,q}(z)$ of $L^{\mathrm{min}}_{f_{p,q}}$ at~$z=0$ equals $z^{p+q-1} \cdot  R_{p,q}(z)$ where $R_{p,q}(0) \neq 0$, hence Yurkevich's conjecture holds.
By looking closely at the indicial polynomials $I_{p,q}(z)$, we made a few additional observations. 
First, the degree of $I_{p,q}(z)$, that is the order of $L^{\mathrm{min}}_{f_{p,q}}$, is equal to $\lfloor (p+q)^2/4\rfloor$ except when $p$ and $q$ are both even and equal, in which case it equals $p^2-1$. 
In particular, the order of $L^{\mathrm{min}}_{f_{p,q}}$ is equal to that of $L^{\mathrm{min}}_{f_{q,p}}$ in the 45 cases. 
{Second, if $p$ is odd or if $p>q$, then $I_{p,q}(z)$ factors
$z^{p+q-1}(z-1)^{p+q-3} \cdots (z-r)^{p+q-2r-1}$},
where $r = \lfloor (p+q)/2\rfloor - 1$; when $p$ and $q$ are
both even and equal, then $I_{p,q}(z)$ equals $z^{2p-1} \cdot  (z-1)^{2p-3} \cdots (z-(p-2))^{3}$.
In the remaining case ($p<q$ with even $p$), the previous pattern may break; for instance,
$I_{2,4}(z) = z^5 (z-1)^3 (z-\frac12)$.
We leave the proof of these experimentally observed formulas, and
their extension to arbitrary $p,q$, as open problems.
(Similar considerations for the exponential generating versions of ${f_{p,q}}$ appear in~\cite[Question~5.2]{BoRiSa24}.)

\appendix
\crefalias{section}{appendix}

\section{Generalized Apéry series}\label{appendix:Apery}

\subsection{A multiple binomial sum}
A direct application of Flajolet's criterion proves the
transcendence of the following
family of power series, extending the one in~\cref{sec:Apery_like}.
\begin{proposition}\label[proposition]{prop:Apery2}
Let  $(p_0,\dots,p_m) \in \N^{m+1}$ with $p_0\ge1$. Then the 
power series 
\[
f_{\mathbf p} (z) =
\sum_{n \geq 0} \left( 
\sum_{k=0}^n 
\binom{n}{k}^{p_0} 
\binom{n+k}{k}^{p_1}
\binom{n+2k}{k}^{p_2}
\dotsm
\binom{n+mk}{k}^{p_m} 
\right) z^n
\]
is transcendental if and only if $p:=p_0+\dots+p_m>2$.
\end{proposition}
\begin{proof}
The case $p>2$ is covered by Flajolet's criterion. Asymptotically, the
$n$th coefficient~$S_ {\mathbf p} (n)$ of~$f_
{\mathbf
p}(z)$ behaves like
\begin{equation}\label{eq:McIntosh}
S_{\mathbf p}(n)\sim \gamma (\pi n)^{
\frac{1-p}2}\beta^n,\quad n\rightarrow\infty
\end{equation}
with $\gamma$ and $\beta$ algebraic~\cite{McIntosh1996}. 
If $p>2$ is odd, the exponent of~$n$ is a negative integer making $f_
{\mathbf p} (z)$ transcendental; if $p$ is even, $\Gamma((3-p)/2)$ is
a rational multiple of $\sqrt{\pi}$ and thus $\gamma\Gamma(
(3-p)/2)\pi^{(1-p)/2}$ is an algebraic multiple of~$\pi^{1-p/2}$,
which is not algebraic for $p>2$.

For $p=1$, the series $f_{(1,0,\dots,0)}$ is the rational~$1/(1-2z)$.

Finally, for $p=2$, either $p_0=2$ and all the other~$p_i$ are~0, in
which case $f_{\mathbf p}=(1-4z)^{-1/2}$, or $p_0=1$, $p_j=1$ for
one~$j\neq0$ and all other $p_i$ are~0. In that case, the power series
can be seen to be algebraic as the diagonal of a bivariate rational
function. Here, the rational
function
\[\frac1{1-z(1+y)(y+(1+y)^j)}
=\sum_{n\ge0}z^n\sum_{k=0}^n\binom{n}{k}y^{n-k}(1+y)^{n+jk},\]
has a Taylor expansion whose coefficient of $z^ny^n$ is the desired
\[
\sum_{k=0}^n\binom{n}{k}\binom{n+jk}{k}.\qedhere
\]
\end{proof}

\subsection{Generic diagonals in more than 2 variables are
transcendental} \label{sec:generic_diags}

	It is well known that the Apéry power series $f_{2,2}$ can be written as the diagonal of a rational power series in 4 variables, e.g., as the diagonal of 
$1/((1-x-y)(1-t-z)-xyzt)$ \cite[Thm.~1.1]{Straub14}
	 or of 
 $1/(((1-z)(1-t)-x)(1-y) - xyzt)$ (see \cref{sec:Apery_like}).
From such a diagonal representation, the approach of analytic
combinatorics in several variables (ACSV) allows to deduce the
asymptotics of the $n$th Apéry number $A_n$ as in~\cref{ex:Apery}, see e.g., \cite[Example~1]{MeSa21}.
More generally, the generalized Apéry power series $f_{p,q}$ can be written as the diagonal of a rational power series in $m=p+q$ variables \emph{with nonnegative coefficients}, e.g., via~\eqref{eq:Apery_as_diag}.
Now, \cite[Result 1]{MeSa21} states that under certain assumptions that hold generically, 
the $n$th coefficient of the diagonal of a rational power series in $m$ variables with nonnegative coefficients grows asymptotically like $\gamma \beta^n n^r$ where 
$r=(1-m)/2$ and $\gamma / (2\pi)^{r} \in \Qbar$.
(This is coherent with the asymptotics~\eqref{eq:McIntosh} in the
proof of \cref{prop:Apery2}.)
Combined with Flajolet's criterion, the formula of \cite[Result 1]
{MeSa21} implies that diagonals of ``generic'' rational functions with
nonnegative coefficients are transcendental if $m>2$; conversely, the
bivariate case is algebraic by the result of P\'olya and
Furstenberg.
Note that the genericity assumption is essential, as shown by the example of the diagonal of 
$1/(1-x-y-yz-xyz)$, equal to $f_{1,1}(z)$ which is algebraic by 
\cref{prop:Apery2}, or by the algebraic examples of \cref{sec:examples-diagonals}.

\section*{Acknowledgements}
Warm thanks go to Jean-Marie Maillard, who provided the final example of~\cref{sec:examples-diagonals}, as well as to
Gilles~Villard, who used his code and Vincent Neiger's code on that
example.
We also thank Gilles Christol and Yves André for many interesting exchanges, in particular regarding \cref{conj:ChristolAndre}.
We are also grateful to Florian Fürnsinn and Sergey Yurkevich, who
gave us constructive advice on a first version of the manuscript. 
This work has been supported by the French--Austrian project 
\href{https://anr.fr/Project-ANR-22-CE91-0007}{EAGLES}
(ANR-22-CE91-0007 \& FWF I6130-N) and by the French project
\href{https://nuscap.gitlabpages.inria.fr/index.html}{NuSCAP}
(ANR-20-CE48-0014).

\bibliographystyle{acm}

\end{document}